\documentclass{amsart}
\usepackage{amssymb}

\newcommand{\Sym}{{\mathfrak{S}}}
\newcommand{\Z}{{\mathbb{Z}}}
\newcommand{\F}{{\mathbb{F}}}
\newcommand{\Q}{{\mathbb{Q}}}
\newcommand{\N}{{\mathbb{N}}}
\newcommand{\R}{{\mathbb{R}}}
\newcommand{\C}{{\mathbb{C}}}
\newcommand{\ba}{{\mathbf{a}}}
\newcommand{\bc}{{\boldsymbol{c}}}
\newcommand{\bD}{{\boldsymbol{D}}}
\newcommand{\cI}{{\mathcal{I}}}
\newcommand{\cO}{{\mathcal{O}}}
\newcommand{\cB}{{\mathcal{B}}}
\newcommand{\cG}{{\mathcal{G}}}
\newcommand{\cL}{{\mathcal{L}}}
\newcommand{\cP}{{\mathcal{P}}}
\newcommand{\cR}{{\mathcal{R}}}
\newcommand{\cS}{{\mathcal{S}}}
\newcommand{\cD}{{\mathcal{D}}}
\newcommand{\cF}{{\mathcal{F}}}
\newcommand{\cC}{{\mathcal{C}}}
\newcommand{\fI}{{\mathfrak{I}}}

\newcommand{\fg}{{\mathfrak{g}}}
\newcommand{\fp}{{\mathfrak{p}}}

\newcommand{\Irr}{{\operatorname{Irr}}}
\newcommand{\Cl}{{\operatorname{Cl}}}
\newcommand{\Res}{{\operatorname{Res}}}
\newcommand{\Ind}{{\operatorname{Ind}}}
\newcommand{\Jind}{{\operatorname{J}}}
\newcommand{\Con}{{\operatorname{Con}}}
\renewcommand{\leq}{\leqslant}
\renewcommand{\geq}{\geqslant}
\renewcommand{\atop}[2]{\genfrac{}{}{0pt}{}{#1}{#2}}

\newtheorem{thm}{Theorem}[section]
\newtheorem{lem}[thm]{Lemma}
\newtheorem{cor}[thm]{Corollary}
\newtheorem{prop}[thm]{Proposition}
\newtheorem{conj}[thm]{Conjecture}

\theoremstyle{definition}
\newtheorem{exmp}[thm]{Example}
\newtheorem{defn}[thm]{Definition}
\newtheorem{abschnitt}[thm]{}

\theoremstyle{remark}
\newtheorem{rem}[thm]{Remark}

\begin{document}

\title{Left cells and constructible representations}

\author{Meinolf Geck}
\address{Institut Girard Desargues, bat. Jean Braconnier, Universit\'e 
Lyon 1, 21 av Claude Bernard, F--69622 Villeurbanne Cedex, France}

\email{geck@igd.univ-lyon1.fr}


\date{April 2004}

\begin{abstract} We consider the partition of a finite Coxeter group $W$ 
into left cells with respect to a weight function $L$. In the equal
parameter case, Lusztig has shown that the representations carried by the 
left cells are precisely the so-called constructible ones. We show that 
this holds for general $L$, assuming that the conjectural properties 
(P1)--(P15) in Lusztig's book on Hecke algebras with unequal parameters 
hold for $W,L$. Our proofs use the idea (Gyoja, Rouquier) that left cell 
representations are projective in the sense of modular representation
theory. This also gives partly new proofs for Lusztig's result in the 
equal parameter case. 
\end{abstract}

\maketitle

\section{Introduction} \label{sec:intro}
Let $W$ be a finite Coxeter group and $L$ be a weight function on $W$,
as in \cite{Lusztig03}. Thus, $L$ is a function $L \colon W\rightarrow\Z$ 
such that $L(ww')=L(w)+L(w')$ whenever $l(ww')= l(w)+l(w')$ where $l$ is 
the length function on $W$. We assume that $L(w)>0$ for all $w\neq 1$. The 
choice of such an $L$ gives rise to a partition of $W$ into left cells; each 
left cell naturally carries a representation of $W$. These representations
play an important role, for example, in the representation theory of 
reductive groups over finite or $p$-adic fields; see \cite{LuBook},
\cite[Chap.~0]{Lusztig03}. In the case where $L=al$ for some $a>0$, 
the representations carried by the left cells are explicitly known: by 
Lusztig \cite{Lusztig86}, they are precisely the constructible 
representations which were defined (and determined) in \cite{Lusztig82b}.

Let us now consider a general weight function $L$ and let us assume that 
the conjectural properties (P1)--(P15) in \cite[Chap.~14]{Lusztig03} hold;
we recall these properties and some of their consequences in Section~2. 
The purpose of this paper is to prove that, in this setting, the left cell 
representations are again the constructible ones, as conjectured by 
Lusztig \cite[22.29]{Lusztig03}. Our arguments also cover the case of
equal parameters, so that we obtain partly new proofs for Lusztig's 
results in \cite{Lusztig86}. 

The main theme of this paper is to use a generalization of a result of 
Rouquier \cite{rou} which shows that left cell representations can be 
interpreted as projective indecomposable representations in the sense of 
modular representation theory. Rouquier's original result was concerned 
with the equal parameter case. The proof of its generalisation to arbitrary 
weight functions relies on (P1)--(P15); see Section~3. In this context,
a crucial role is played by a general result from \cite{GeRo} which shows 
that the columns of the ``decomposition matrix'' are linearly independent; 
see (\ref{remcent}). This result yields, for example, the fact that the 
representations carried by the various left cells of $W$ are linearly 
independent in the appropriate Grothendieck group; see Remark~\ref{cc1}.

See Gyoja \cite{Gy} where the connection between left cells and modular 
representations first appeared. 

Further general results on left cells, constructible representations and 
families will be discussed in Sections~4 and~5. 

In order to deal with type $B_m$ (Section~6), we need a purely
combinatorial identity about the constructible representations; see 
Proposition~\ref{linind}. We shall give here an elementary proof of that 
identity which uses the fact that the constructible representations
can be identified with the canonical basis of a certain irreducible 
representation of the Lie algebra $\mathfrak{sl}_{n+1}(\C)$, as shown by
Leclerc--Miyachi \cite{HyLe}. (I wish to thank Bernard Leclerc for his help
with the proof of Proposition~\ref{linind}.)

Our methods provide new proofs for type $D_m$ and $B_m$ with equal 
parameters.

Finally, in Section~7, we discuss the computation of left cell 
representations in groups of exceptional type. For groups of type 
$E_6$, $E_7$, $E_8$ and $F_4$ with equal parameters, we show that
Lusztig's proof in \cite{Lusztig86}, which partly relies on deep 
results about representations of reductive groups over finite fields 
\cite{LuBook}, can be replaced by arguments which only rely on (P1)--(P15) 
and some explicit computations with the character tables in \cite{ourbuch}.

\section{Lusztig's conjectures} \label{conj}

Let $(W,S)$ be a Coxeter system where $W$ is finite; let $L\colon W
\rightarrow \Z$ be a weight function such that $L(s)>0$ for all $s\in S$.
Let $A={\Z}[v,v^{-1}]$ where $v$ is an indeterminate. We set $v_w=v^{L(w)}$ 
for all $w\in W$. Let $H$ be the generic Iwahori--Hecke algebra corresponding 
to $(W,S)$ with parameters $\{v_s \mid s \in S\}$. Thus, $H$ has an $A$-basis
$\{T_w\mid w \in W\}$ and the multiplication is given by the rule
\[ T_sT_w=\left\{\begin{array}{cl} T_{sw} &\quad\mbox{if $l(sw)=l(w)+1$},\\
T_{sw}+(v_s-v_s^{-1})T_w & \quad \mbox{if $l(sw)=l(w)-1$},\end{array}
\right.\]
where $s\in S$ and $w\in W$. 
Let $\{c_w\mid w\in W\}$ be the ``new'' basis of $H$ defined in 
\cite[Theorem~5.2]{Lusztig03}. We have $c_w=T_w+\sum_{y} p_{y,w}\, T_y$
where $p_{y,w} \in v^{-1}{\Z}[v^{-1}]$ and $p_{yw}=0$ unless $y<w$
in the Bruhat--Chevalley order. Given $x,y\in W$, we can write 
\[ c_xc_y=\sum_{z\in W} h_{x,y,z}c_z \qquad \mbox{where $h_{x,y,z}\in A$}.\]
We usually work with the elements $c_w^\dagger$
obtained by applying the unique $A$-algebra involution $H \rightarrow H$, 
$h \mapsto h^\dagger$ such that $T_s^\dagger=-T_s^{-1}$ for any $s\in S$; 
see \cite[3.5]{Lusztig03}. 

We refer to \cite[Chap.~8]{Lusztig03} for the definition of the
preorders $\leq_{\cL}$, $\leq_{\cR}$, $\leq_{\cL\cR}$ and the corresponding
equivalence relations $\sim_{\cL}$, $\sim_{\cR}$, $\sim_{\cL\cR}$. (See also
Remark~\ref{rem0} below.) The equivalence classes with respect to these 
relations are called left, right and two-sided cells, respectively.

Let $\Gamma$ be a left cell of $W$. As in \cite[21.1]{Lusztig03}, the 
$A$-submodule 
\[ [\Gamma]_A:=\sum_{y\in \Gamma} Ac_y^\dagger \]
of $H$ can be regarded as an $H$-module by the rule
$c_x^\dagger \cdot c_y^\dagger=\sum_{z\in \Gamma} h_{x,y,z}c_z^\dagger$
where $x\in W$, $y\in \Gamma$.
By change of scalars ($v\mapsto 1$) this gives rise to an ${\C}[W]$-module
which we simply denote by $[\Gamma]$. The representation $[\Gamma]$ is called
a {\em left cell representation} of $W$. 

\begin{conj}[Lusztig \protect{\cite[22.29]{Lusztig03}}] 
\label{conjC} Let $\Gamma$ be a left cell of $W$. Then $[\Gamma]$ is a 
constructible representation in the sense of \cite[22.1]{Lusztig03}.
Furthermore, all constructible representations arise in this way.
\end{conj}

We recall the definition of constructible representations in (\ref{fam}).
Conjecture~\ref{conjC} is already  known to be true if $L=al$ for some 
$a>0$, see \cite{Lusztig86} (for $W$ a finite Weyl group) and \cite{Al} 
(for $W$ of type $H_4$). 

We will show in Sections~6 and~7 that Conjecture~\ref{conjC} is true, if 
the general conjectural properties (P1)--(P15) formulated by Lusztig in 
\cite[Chap.~14]{Lusztig03} hold. This also yields a partly new proof for 
finite Weyl groups and $L=al$ for some $a>0$ (where (P1)--(P15) are known 
to hold; see \cite[Chap.~15]{Lusztig03}.) Assuming (P1)--(P15), we already 
know that every constructible representation is of the form $[\Gamma]$ for 
some left cell $\Gamma$; see \cite[Lemma~22.2]{Lusztig03}. So it remains 
to show that $[\Gamma]$ is constructible for every left cell $\Gamma$.

To state (P1)--(P15), we have to introduce some further notation.

For $z\in W$, there is a unique integer $\ba(z)\geq 0$ such that 
\begin{alignat*}{2}
h_{x,y,z} &\in v^{\ba(z)}{\Z}[v^{-1}] &&\qquad \mbox{for all $x,y\in W$},\\
h_{x,y,z} &\not\in v^{\ba(z)-1}{\Z}[v^{-1}] &&\qquad \mbox{for some $x,y\in 
W$};
\end{alignat*}
see \cite[13.6]{Lusztig03}. Given $x,y,z\in W$, we define 
$\gamma_{x,y,z^{-1}}\in \Z$ by the condition  
\[ h_{x,y,z}=\gamma_{x,y,z^{-1}}v^{\ba(z)}+
\mbox{strictly smaller powers of $v$}.\]
Finally, for $z\in W$, we define an integer $\Delta(z)\geq 0$ and a non-zero
integer $n_z\in \Z$ by the condition
\[ p_{1,z}=n_zv^{-\Delta(z)}+\mbox{strictly smaller  powers of $v$}; 
\quad \mbox{see \cite[14.1]{Lusztig03}}.\]

\begin{conj}[Lusztig \protect{\cite[14.2]{Lusztig03}}] \label{Pconj} 
Let $\cD=\{z\in W \mid \ba(z)=\Delta(z)\}$.  Then the  following properties 
hold.
\begin{itemize}
\item[\bf P1.] For any $z\in W$ we have $\ba(z)\leq \Delta(z)$.
\item[\bf P2.] If $d \in \cD$ and $x,y\in W$ satisfy $\gamma_{x,y,d}\neq 0$,
then $x=y^{-1}$.
\item[\bf P3.] If $y\in W$, there exists a unique $d\in \cD$ such that 
$\gamma_{y^{-1},y,d}\neq 0$.
\item[\bf P4.] If $z'\leq_{\cL\cR} z$ then $\ba(z')\geq \ba(z)$. Hence, if
$z'\sim_{\cL\cR} z$, then $\ba(z)=\ba(z')$.
\item[\bf P5.] If $d\in \cD$, $y\in W$, $\gamma_{y^{-1},y,d}\neq 0$, then
$\gamma_{y^{-1},y,d}=n_d=\pm 1$.
\item[\bf P6.] If $d\in \cD$, then $d^2=1$.
\item[\bf P7.] For any $x,y,z\in W$, we have $\gamma_{x,y,z}=\gamma_{y,z,x}$.
\item[\bf P8.] Let $x,y,z\in W$ be such that $\gamma_{x,y,z}\neq 0$. Then 
$x\sim_{\cL}y^{-1}$, $y \sim_{\cL}z^{-1}$, $z\sim_{\cL} x^{-1}$.
\item[\bf P9.] If $z'\leq_{\cL} z$ and $\ba(z')=\ba(z)$, then $z'\sim_{\cL}z$.
\item[\bf P10.] If $z'\leq_{\cR} z$ and $\ba(z')=\ba(z)$, then $z'\sim_{\cR}z$.
\item[\bf P11.] If $z'\leq_{\cL\cR} z$ and $\ba(z')=\ba(z)$, then
$z'\sim_{\cL\cR}z$.
\item[\bf P12.] Let $I\subset S$. If $y\in W_I$, then $\ba(y)$ computed
in terms of $W_I$ is equal to $\ba(y)$ computed in terms of $W$.
\item[\bf P13.] Any left cell $\Gamma$ of $W$ contains a unique element
$d\in \cD$. We have $\gamma_{x^{-1},x,d}\neq 0$ for all $x\in \Gamma$.
\item[\bf P14.] For any $z\in W$, we have $z \sim_{\cL\cR} z^{-1}$.
\item[\bf P15.] Let $v'$ be a second indeterminate and let $h_{x,y,z}'\in
{\Z}[v',v'^{-1}]$ be obtained from $h_{x,y,z}$ by the substitution $v
\mapsto v'$. If $x,x',y,w\in W$ satisfy $\ba(w)=\ba(y)$, then 
\[\sum_{y' \in W} h_{w,x',y'}'\,h_{x,y',y}=\sum_{y'\in W} h_{x,w,y'}\,
h_{y',x',y}'.\]
\end{itemize}
\end{conj}

By \cite[Chap.~15]{Lusztig03}, the above conjectures hold if $L=al$ for
some $a>0$. We shall assume from now on that the above conjectures hold
for $W,L$. Then we can define a new algebra $J$ over $\Z$ as in
\cite[Chap.~18]{Lusztig03}. As a ${\Z}$-module, $J$ is free with
a basis $\{t_w\mid w\in W\}$. The multiplication is defined by 
\[ t_xt_y=\sum_{z\in W} \gamma_{x,y,z}\, t_{z^{-1}}.\]
This multiplication is associative and we have an identity element given
by $1_J=\sum_{d\in \cD} n_dt_d$. Furthermore, we have a homomorphism
of $A$-algebras $\phi \colon H \rightarrow J_A=A\otimes_{\Z} J$ given by
\[ \phi(c_w^\dagger)=\sum_{\atop{z\in W,d\in \cD}{\ba(z)=\ba(d)}}
h_{w,d,z} \hat{n}_z\, t_z \qquad (w\in W),\]
where $\hat{n}_z$ is defined as follows. For any $z\in W$, we set 
$\hat{n}_z=n_d$ where $d$ is the unique element of $\cD$ such that 
$d\sim_{\cL} z^{-1}$ and $n_d=\pm 1$; see (P5). 

\begin{rem} \label{rem1} Let $w,z\in W$ and assume that there exists
some $d\in \cD$ such that $\ba(z)=\ba(d)$ and $h_{w,d,z} \neq 0$. Then we
have $z \leq_{\cL} d$. So (P9) implies that $z\sim_{\cL} d$.
But then (P13) shows that $d$ is the unique element of $\cD$ in the same
left cell as $z$. Thus, we have in fact 
\[ \phi(c_w^\dagger)=\sum_{z\in W} \hat{n}_{z}\,h_{w,d_z,z}\, t_z\]
where, for any $z\in W$, we denote by $d_z$ the unique element of
$\cD$ in the same left cell as $z$. We have 
\[\pm v^e\det(\phi) \in 1+v{\Z}[v] \qquad \mbox{where $e=\sum_w \ba(w)$}.\]
Indeed, let $w\in W$. Then, by (P2), (P4), (P5) and (P12), 
we have 
\begin{align*} v^{\ba(w)}\,\phi(c_w^\dagger) =\pm t_w&+
\mbox{$v{\Z}[v]$-combination of elements $t_z$ where $\ba(z)=\ba(w)$}\\
&+\mbox{$A$-combination of elements $t_x$ where $\ba(x)>\ba(w)$}.
\end{align*}
Hence we see that the matrix of $\phi$ has a block triangular shape,
when we order the elements of $W$ according to increasing value of $\ba$.
Furthermore, inside a diagonal block, the coefficients have leading term 
$\pm v^{\ba(w)}$ on the diagonal and strictly smaller leading term off 
the diagonal.
\end{rem}

\begin{rem} \label{rem0} 
Using $J$, the relations $\sim_{\cL}$, $\sim_{\cR}$, $\sim_{\cL\cR}$ can be 
characterized as follows (see \cite[Prop.~18.4]{Lusztig03}):
\begin{itemize}
\item We have $x\sim_{\cL} y$ if and only if $t_{x}t_{y^{-1}}\neq 0$. 
\item We have $x\sim_{\cR} y$ if and only if $t_{x^{-1}}t_{y}\neq 0$. 
\item We have $x\sim_{\cL\cR} y$ if and only if $t_xt_wt_{y}\neq 0$ for
some $w\in W$. 
\end{itemize}
\end{rem}

\section{Cells and idempotents} \label{blocks}

The purpose of this section is to generalize the main result of Rouquier 
\cite{rou} to the unequal parameter case, assuming that (P1)--(P15) hold.

Let $B_0$ be the set of two-sided cells in $W$. We set 
\[ t_{\bc}=\sum_{d\in \cD \cap \bc} n_dt_d\qquad\mbox{for $\bc \in B_0$}.\]

\begin{lem} \label{lem31} We have $1_J=\sum_{\bc \in B_0} t_{\bc}$. 
Furthermore, $t_{\bc}t_{\bc'}=\delta_{\bc \bc'} t_{\bc}$ for any $\bc,\bc'
\in B_0$.  Thus, $\{t_{\bc} \mid \bc \in B_0\}$ is a set of mutually
orthogonal, central idempotents in $J$.
\end{lem}

\begin{proof} For $\bc \in B_0$, let $J_{\bc}=\langle t_w \mid w\in \bc
\rangle_{\Z}$. By (P8), $J_{\bc}$ is a two-sided ideal in $J$; we have
$J=\bigoplus_{\bc \in B_0} J_{\bc}$. This yields a unique decomposition 
$1_J=\sum_{\bc \in B_0} e_{\bc}$ where $e_{\bc} \in J_{\bc}$. Here, 
$\{e_{\bc} \mid \bc \in B_0\}$ is a set of mutually orthogonal, central 
idempotents. For $\bc \in B_0$, let us write $e_{\bc}=\sum_{w\in \bc} 
a_{\bc,w}\, t_w$ where $a_{\bc,w} \in \Z$. Since $1_J=\sum_{d\in \cD} 
n_dt_d$, we conclude that $a_{\bc,w}=n_w$ for $w\in \cD \cap \bc$, and 
$0$ otherwise. Thus, we have $e_{\bc}=t_{\bc}$ as required.
\end{proof}

\begin{lem} \label{lem30} Let $\Gamma$ be a left cell and let $\cD\cap 
\Gamma=\{d\}$. Then we have $n_dt_wt_d=t_w$ for any $w\in \Gamma$.
Furthermore, we have   
\[ Jt_d=\langle t_y \mid y \in \Gamma\rangle_{\Z} \qquad \mbox{and}\qquad
(n_dt_d)^2=n_dt_d.\]
\end{lem}

\begin{proof} Let $y\in W$. We have $t_yt_d=\sum_{x\in W} \gamma_{y,d,x} 
t_{x^{-1}}$ where $\gamma_{y,d,x}\in \Z$. If $x=y^{-1}$ and $y\in \Gamma$, 
then $\gamma_{y,d,x}=\gamma_{x,y,d}=n_d=\pm 1$ by (P5), (P7), (P12). Now 
consider any $x\in W$ and assume that $\gamma_{y,d,x}\neq 0$. We must show 
that $x=y^{-1}$ and $y\in \Gamma$. Now, by (P7), we have $\gamma_{x,y,d}=
\gamma_{y,d,x} \neq 0$. By (P2) and (P8), this implies $x=y^{-1}$ and 
$y\in \Gamma$, as required.
\end{proof}

We shall see that each $n_dt_d$ actually is a primitive idempotent in $J$. In 
fact, this will even work when we extend scalars from $\Z$ to suitable larger 
rings. In order to describe the required conditions on such a larger ring,
we recall the following constructions.

\begin{abschnitt} {\bf The simple $J_{\C}$-modules.} \label{simpJ}
Upon substituting $v\mapsto 1$, the algebra $H$ specialises to ${\Z}[W]$. 
Hence, extending scalars from $\Z$ to $\C$, we obtain an isomorphism of 
$\C$-algebras
\[ \phi_1\colon {\C}[W]  \rightarrow J_{\C}, \qquad \mbox{where $J_{\C}=
{\C} \otimes_{\Z} J$};\]
see \cite[20.1]{Lusztig03}. Since ${\C}[W]$ is split semisimple,  we can 
now conclude that $J_{\C}$ also is split semisimple. 

Let $\Irr(W)$ be the set of simple ${\C}[W]$-modules up to isomorphism.
For any ${\C}[W]$-module $E$, we denote the corresponding $J_{\C}$-module by
$E_{\spadesuit}$. Thus, $E_{\spadesuit}$ coincides with $E$ as an
$\C$-vectorspace and the action of $a\in J_{\R}$ on $E_{\spadesuit}$ is 
the same as the action of $\phi_1^{-1}(a)$ on $E$; see \cite[20.2]{Lusztig03}. 
Then we have 
\[ \Irr(J_{\C})=\{E_{\spadesuit} \mid E \in \Irr(W)\},\]
where $\Irr(J_{\C})$ is the set of simple $J_{\C}$-modules up to isomorphism.
We shall also need the fact that $J_{\C}$ is a symmetric algebra, with trace
form $\tau \colon J_{\C} \rightarrow \C$ given by $\tau(t_z)=n_z$ if
$z\in \cD$ and $\tau(t_z)=0$ otherwise. We have $\tau(t_xt_y)=
\delta_{xy,1}$ for any $x,y\in W$; see \cite[20.1]{Lusztig03}. By general 
results on split semisimple symmetric algebras (see 
\cite[Chap.~7]{ourbuch}), we have
\[ \tau(t_w)=\sum_{E\in \Irr(W)} \frac{1}{f_E}\, \mbox{tr}(t_w,E_\spadesuit)
\qquad \mbox{for all $w\in W$},\]
where $0\neq f_E \in \C$ for all $E \in \Irr(W)$. By 
\cite[Lemma~20.13]{Lusztig03}, we have in fact 
\[ f_E \in \R \qquad \mbox{and}\qquad f_E>0.\]
\end{abschnitt}

\begin{abschnitt} {\bf The simple $H_K$-modules.} \label{simpH}
Let us extend scalars from $A$ to $K={\C}(v)$. Then we obtain an 
isomorphism  of $K$-algebras 
\[ \phi_K\colon H_K\rightarrow J_K \qquad \mbox{where $H_K=K\otimes_A H$
and $J_K=K \otimes_{\Z} J$};\]
see \cite[20.1]{Lusztig03}. Given a ${\C}[W]$-module $E$, the 
$J_{\C}$-module structure on $E_\spadesuit$ extends in a natural way to a 
$J_K$-module structure on $E_v:=K\otimes_{\C} E_\spadesuit$.  Then we can 
also regard $E_v$ as an $H_K$-module via $\phi_K$. We have 
\[\Irr(H_K)=\{E_v \mid E \in \Irr(W)\},\]
where $\Irr(H_K)$ is the set of simple $H_K$-modules up to isomorphism.
We have $\mbox{tr}(T_w,E_v) \in {\C}[v,v^{-1}]$ for all $w\in W$ and 
every ${\C}[W]$-module $E$. This yields the following direct relation 
between $E$ and $E_v$ (see \cite[20.3]{Lusztig03}):
\[ \mbox{tr}(w,E)=\mbox{tr}(T_w,E_v)\mid_{v=1} \qquad \mbox{for all 
$w\in W$}.\]
The above direct relation between $\Irr(W)$ and $\Irr(H_K)$ can also
be established without reference to the algebra $J$ and (P1)--(P15); see
\cite[8.1.7 and 9.3.5]{ourbuch}. 

Let $E\in \Irr(W)$. As in \cite[Prop.~20.6]{Lusztig03}, define an integer
$\ba_E \geq 0$ by the condition
\begin{alignat*}{2}
v^{{\ba}_E}\mbox{tr}(T_w,E_v) &\in {\C}[v] &&\qquad \mbox{for all $w\in W$},
\\ v^{{\ba}_E-1}\mbox{tr}(T_w,E_v) &\not\in {\C}[v] &&\qquad \mbox{for some 
$w\in W$}.
\end{alignat*}
Note that $\ba_E$ depends on the choice of $L$. Now let 
\[c_E:= \frac{1}{\dim E}\sum_{w\in W} \mbox{tr}(T_w,E_v)\,
\mbox{tr}(T_{w^{-1}},E_v) \in {\C}[v,v^{-1}].\]
In \cite{ourbuch}, this is called the {\em Schur element} associated to $E$. 
Then, by \cite[Cor.~20.11]{Lusztig03}, we have 
\[ c_E=f_E v^{-2\ba_E}+\mbox{combination of higher powers of $v$},\]
where $f_E$ is the positive real number introduced in (\ref{simpJ}).
\end{abschnitt}

\begin{defn} \label{cadapt} By Lemma~\ref{lem31}, we have a partition
\[ \Irr(W)=\coprod_{\bc \in B_0} \Irr(W,\bc)\]
where $\Irr(W,\bc)=\{E \in \Irr(W) \mid t_{\bc}E_{\spadesuit}\neq 0\}$.
The sets $\Irr(W,\bc)$ may be called the ``blocks'' of $\Irr(W)$. Now 
fix $\bc \in B_0$ and let $R\subseteq {\C}(v)$ be a notherian subring.
We say that $R$ is {\em $\bc$-adapted} if the following condition holds:
\begin{equation*}
R\cap \Big\{\sum_{E\in \Irr(W,\bc)} 
\frac{n_E}{f_E} \;\Big|\; n_E \in \Z\Big\} \subseteq \Z.\tag{$*_{\bc}$} 
\end{equation*}
\end{defn}

\begin{exmp} \label{exp1} Here are some examples of $\bc$-adapted 
rings $R$.

(a) The subring $\Z\subseteq {\C}(v)$ clearly satisfies ($*_\bc$) for
any $\bc \in B_0$.

(b) Following Rouquier \cite{rou}, let $\fI_+=\{1+vf\mid f \in 
{\Z}[v]\}$ and $\cO=\{f/g \in {\C}(v) \mid  f\in A,g\in \fI_+\}$. 
Then $\cO$ is a principal ideal domain; see \cite[p.~1040]{rou}. 
It is easily checked that $\cO$ satisfies ($*_\bc$) for any $\bc \in B_0$.

(c) Let $\bc \in B_0$ and assume that there is a prime number $p$ such that 
$f_E$ is a power of $p$ for any $E\in \Irr(W,\bc)$. Let $A_p=\{f/g \in 
{\C}(v) \mid  f,g\in A, g\not\in pA\}$. Then $A_p$ is a noetherian local 
ring whose maximal ideal is the principal ideal generated by $p$. It is 
readily checked that ($*_\bc$) holds for $A_p$. Rings of this type have been 
used by Gyoja \cite{Gy} in the study of left cells in the equal parameter 
case.
\end{exmp}

\begin{thm}[See Rouquier \cite{rou} in the equal parameter case] 
\label{prop32} Let $\bc \in B_0$ and $R\subseteq {\C}(v)$ be a $\bc$-adapted
subring. By extension of scalars, we obtain an $R$-algebra $J_R=R 
\otimes_{\Z} J$. Then the following hold.
\begin{itemize}
\item[(a)] Let $d\in \cD\cap\bc$. Then $n_dt_d$ is a primitive idempotent 
in $J_{R}$.
\item[(b)] $t_{\bc}$ is a primitive idempotent in the center of $J_{R}$.
\end{itemize}
\end{thm}

\begin{proof} Let $K={\C}(v)$ and consider $J_{R}$ as a subalgebra 
of $J_K=K \otimes_{\Z} J$. Now, for any $E \in \Irr(W)$, we can extend the 
$J_{\C}$-module structure on $E_\spadesuit$ in a natural way to a 
$J_K$-module structure on $K\otimes_{\C} E_\spadesuit$. The trace form
$\tau$ also extends to a trace form on $J_K$ which we denote by the same
symbol. Thus, we have 
\[ \tau(h)=\sum_{E \in \Irr(W)} \frac{1}{f_E}\, \mbox{tr}(h,K\otimes_{\C}
E_\spadesuit) \qquad \mbox{for any $h\in J_K$}.\]

(a) We begin with the following general remark. Let $e\in J_{R}$ be any 
primitive idempotent such that $e=et_{\bc}$. We regard $e$ as an element
of $J_K$. Now, for $E\in \Irr(W)$, we have $eE_{\spadesuit}=0$ unless $E
\in \Irr(W,\bc)$. Furthermore, the term $\mbox{tr}(e,K\otimes_{\C} 
E_\spadesuit)$ is a non-negative integer which is less than or equal to 
$\dim E$. Since $e=et_{\bc}$, there is at least one $E\in\Irr(W,\bc)$ such 
that $\mbox{tr}(e,K\otimes_{\C} E_\spadesuit) >0$. Now recall from the
discussion in (\ref{simpJ}) that $f_E$ is a positive real number for all~$E
\in \Irr(W)$. Combining all these pieces of information, the above formula 
shows that $\tau(e)$ is a positive real number. 

On the other hand, the defining formula for $\tau$ and the fact that $e$ 
lies in $J_{R}$ show that $\tau(e)\in R$. Hence, using ($*_\bc$), we conclude
that 
\[ \tau(e)>0 \quad \mbox{and} \quad \tau(e) \in R\cap 
\Big\{\sum_{E\in \Irr(W,\bc)} \frac{n_E}{f_E} \;\Big|\; n_E \in \Z\Big\}
\subseteq \Z.\]
In particular, this proves that $\tau(e)\geq 1$ for any primitive idempotent
$e\in J_R$ such that $e=et_{\bc}$. Now let $d \in \cD$. By 
Lemma~\ref{lem30}, $n_dt_d$ is an idempotent; furthermore, we have 
$\tau(n_dt_d)=1$. Let us write $n_dt_d=e_1+\cdots +e_n$ where $e_i$
are orthogonal primitive idempotents in $J_R$. (This is possible since
$R$ is noetherian.) Since $d\in \bc$, we have $e_it_{\bc}=e_i$ for all~$i$.
The above discussion shows that $1=\tau(n_dt_d)=\tau(e_1)+\cdots +\tau(e_n) 
\geq n$. Thus, we must have $n=1$ and so $n_dt_d$ is primitive.

(b) Let $\bc \in B_0$ and take some $d \in \cD \cap \bc$. We consider 
a decomposition $1_J=e_1+\cdots + e_r$, where the $e_i$ are mutually 
orthogonal, central primitive idempotents in $J_{R}$. (This exists since
$R$ is noetherian.) We have $n_dt_d=n_dt_d1_J=\sum_i n_dt_de_i$. Since 
$n_dt_d$ is a primitive idempotent, there exists a unique $j$ such that
$n_dt_de_j=n_dt_d$ and $n_dt_de_i=0$ for $i\neq j$.  Thus, we have 
\[t_d\in J_{R}e_j \qquad \mbox{where} \qquad t_de_j\neq 0.\]
Now consider any element $d'\in \cD\cap \bc$. The fact that $d,d'$ lie
in the same two-sided cell means that there exists some $w\in W$ such that 
$t_dt_wt_{d'}\neq 0$; see Remark~\ref{rem0}. Since $t_dt_wt_{d'}
\in J_{R}e_j$, we also have $t_dt_wt_{d'}e_j\neq 0$ and so $t_{d'}e_j
\neq 0$. By the previous argument, this implies $t_{d'}\in J_{R}e_j$. 

Thus, we have shown that $t_{d'}\in J_{R}e_j$ for all $d' \in \cD\cap 
\bc$ and so $t_{\bc} \in J_{R}e_j$. Since $e_j$ is a primitive 
idempotent, we conclude that $t_{\bc}=e_j$, as required.
\end{proof}

Recall that, given a left cell $\Gamma$ of $W$, we have a corresponding
$H$-module $[\Gamma]_A=\sum_{y\in \Gamma} Ac_y^\dagger$. If $R$ is any 
(commutative) ring containing $A$, let $H_R=R\otimes_A H$. Then the 
$H$-module structure on $[\Gamma]_A$ extends in a natural way to an 
$H_R$-module structure on $[\Gamma]_R:=\sum_{y\in \Gamma} Rc_y^\dagger$. 

\begin{cor} \label{lem34} Let $R \subseteq {\C}(v)$ be a noetherian 
subring such that $\cO\subseteq R$. Let $\Gamma$ be a left cell of $W$. 
Then $[\Gamma]_R$ is a projective $H_R$-module. Furthermore, if 
$\Gamma$ is contained in $\bc \in B_0$ and $R$ is $\bc$-adapted, then
$[\Gamma]_R$ is indecomposable.
\end{cor}

\begin{proof} By extension of scalars, we obtain a homomorphism 
$\phi_R \colon H_R \rightarrow J_R$ of $R$-algebras. By Remark~\ref{rem1}, 
$\det(\phi)$ is invertible in $\cO \subseteq R$; hence $\phi_R$ is an 
isomorphism. Now let $\cD \cap \Gamma=\{d\}$. By Lemma~\ref{lem30}, we have 
\[ J_Rt_d=\langle t_y \mid y\in \Gamma\rangle_R.\]
Now, since $n_dt_d$ is an idempotent, $J_Rt_d$ is a projective 
$J_R$-module. Using $\phi_R$, we may also regard $J_Rt_d$ as an 
$H_R$-module. Thus, the action of $h \in H_R$ on $J_Rt_d$ is given by 
$h *t_y:=\phi(h)t_y$ ($y \in \Gamma$).  Since $\phi_R$ is an isomorphism, 
the resulting $H_R$-module is projective. Finally, consider the $R$-linear 
bijection 
\[ \theta \colon J_Rt_d \rightarrow [\Gamma]_R, \qquad t_y \mapsto 
\hat{n}_y\, c_y^\dagger.\]
By a computation analogous to that in \cite[18.10]{Lusztig03}, we obtain 
\[c_w^\dagger * t_y=\phi(c_w^\dagger)t_y=\hat{n}_y \sum_{u\in \Gamma} 
h_{w,y,u}\,\hat{n}_u t_u \qquad \mbox{for any $w\in W$, $y\in \Gamma$}.\]
Applying $\theta$ yields $\theta(c_w^\dagger *t_y)=c_w^\dagger 
\cdot \theta(t_y)$. Thus, $\theta$ is an $H_R$-module isomorphism.
Consequently, $[\Gamma]_R$ is a projective $H_R$-module.

If $\Gamma\subseteq \bc$ and $R$ is $\bc$-adapted, then $n_dt_d$ is a 
primitive idempotent by Theorem~\ref{prop32} and so $[\Gamma]_R$ is seen 
to be indecomposable.
\end{proof}

Note that the statement of Corollary~\ref{lem34} only involves the notion
of a left cell of $W$. The ring $J$ and the properties (P1)--(P15) are 
needed in the proof. It would be very interesting to find a more elementary
proof.

\begin{abschnitt} {\bf Left cells and decomposition numbers.}
\label{remcent}  Assume that we have a discrete valuation ring 
$R \subseteq {\C}(v)$ such that $\cO \subseteq R$, where $\cO$ is the ring 
in Example~\ref{exp1}(b). Assume, furthermore, that 
\begin{equation*}
H_F=F\otimes_R H_R \mbox{ is split semisimple and $H_k=k \otimes_R H_R$ 
is split},\tag{0}
\end{equation*}
where $F$ is the field of fractions of $R$ and $k$ is the residue field 
of $R$. As pointed out in \cite[\S 3]{GeRo2}, these assumptions imply that 
the  Krull--Schmidt Theorem holds for $H_R$-modules which are finitely 
generated and free over $R$; furthermore, idempotents can be lifted
from $H_k$ to $H_R$. (That is, we don't have to pass to a completion of 
$R$ as is usually done in the modular representation theory of finite 
groups and associative algebras.) Now we are in the general setting of 
\cite[\S 7.5]{ourbuch}. The canonical map $R \rightarrow k$ induces a 
decomposition map
\[ d_R\colon R_0(H_F)\rightarrow R_0(H_k)\]
between the Grothendieck groups of finite-dimensional representations
of $H_F$ and $H_k$, respectively. Let $D_R$ be the corresponding 
decomposition matrix. Thus, $D_{R}$ has rows labelled by $\Irr(H_F)$ 
and columns labelled by $\Irr(H_k)$. Note that the extension of scalars
from $F$ to $K={\C}(v)$ induces a natural bijection $\Irr(H_F)
\stackrel{\sim}{\rightarrow} \Irr(H_K)$ (since $H_F$ is assumed to
be split semisimple). Thus, as in (\ref{simpH}), we also have 
a bijection $\Irr(H_F) \stackrel{\sim}{\rightarrow} \Irr(W)$. 

Now the entries of $D_R$ are given as follows. Consider the projective 
indecomposable $H_R$-modules (PIM's for short). Every PIM has a unique 
simple quotient, which is a simple $H_k$-module. In fact, associating to 
each PIM its simple quotient defines a bijection between isomorphism classes 
of PIM's and $\Irr(H_k)$. For each $V\in \Irr(H_k)$, choose a PIM $P_V$ with 
simple quotient $V$. Thus, $\{P_V\mid V \in \Irr(H_k)\}$ is the set of all 
PIM's of $H_R$, up to isomorphism. By Brauer reciprocity (see 
\cite[Theorem~7.5.2]{ourbuch}), the coefficients in a fixed column of 
$D_{R}$ give the expansion of the 
corresponding PIM (viewed as an $H_F$-module by extension of scalars 
from $R$ to $F$) in terms of the irreducible ones. Thus, we have
\begin{equation*}
D_{R}=\Big([E_v:P_V]\Big)_{E\in \Irr(W),V\in \Irr(H_k)} \tag{1}
\end{equation*}
where $\Irr(W)$ indexes the rows and $\Irr(H_k)$ indexes the columns.
Here, we denote by $[E_v :P_V]$ the multiplicity of $E_v\in \Irr(H_K)$ 
as a simple component of $P_V$ (viewed as an $H_K$-module by extension of 
scalars). Now, since $H$ is a symmetric algebra with a ``reduction-stable'' 
center in the sense of \cite[Def.~7.5.5]{ourbuch}, we know by a general 
argument (due to Rouquier and the author, see \cite[Theorem~7.5.6]{ourbuch})
that 
\begin{equation*}
\mbox{the columns of $D_R$ are linearly independant over $\Q$}.\tag{2}
\end{equation*}
Now assume, moreover, that $R$ is $\bc$-adapted for every two-sided cell 
$\bc$. By Corollary~\ref{lem34}, $[\Gamma]_R$ is a PIM. Since $W$ is the 
union of all left cells, we obtain all PIM's in this way. For each 
$V\in \Irr(H_k)$, choose a left cell $\Gamma_V$ such that $[\Gamma_V]_R$ 
has simple quotient $V$. Thus, $\{[\Gamma_V]_R \mid V \in \Irr(H_k)\}$ is 
the set of all PIM's of $H_R$, up to isomorphism. Thus, we finally obtain
\begin{equation*}
D_{R}=\Big([E:[\Gamma_V]]\Big)_{E\in \Irr(W),V\in \Irr(H_k)}. \tag{3}
\end{equation*}
Here, we denote by $[E :[\Gamma]]$ the multiplicity of $E$ as a simple
component of $[\Gamma]$; note that $[E_v:[\Gamma]_K]]=[E:[\Gamma]]$ for 
any $E\in \Irr(W)$ and any left cell~$\Gamma$ (see, for example, the
argument in \cite[20.5]{Lusztig03}).

The above discussion also shows that, given any 
left cell $\Gamma$, there exists a unique $V \in \Irr(H_k)$ such that
$[\Gamma]_R \cong [\Gamma_V]_R$. Thus, we also have $[\Gamma]\cong 
[\Gamma_V]$ as ${\C}[W]$-modules.  
\end{abschnitt}

The above idea of relating left cells and decomposition matrices has 
already been considered by Gyoja \cite{Gy} in the case of equal
parameters. 

\begin{abschnitt} {\bf Construction of $(F,R,k)$}. \label{split} We now 
give a general construction of a discrete valuation ring $R\supseteq \cO$ 
satisfying the hypotheses in (\ref{remcent})(0), where $k$ has 
characteristic $p>0$. 

Let $F_0$ be a finite extension of $\Q$; let $R_0$ be the ring of algebraic
integers in $F_0$. Let $\fp$ be a prime ideal such that $p\in \fp$. Then 
$k=R_0/\fp$ is a finite field of characteristic~$p$. Now we may choose 
$F_0$ large enough so that $k_0 \otimes_{\Z} J$ is split. (This is clearly 
possible.) Then $J_k=k\otimes_{\Z} J$ also is split where $k=k_0(v)$. Now 
$\phi_{k} \colon H_k\rightarrow J_k$ is an isomorphism (see 
Remark~\ref{rem1}). So we conclude that $H_k$ is split. On the other hand, 
$H_F$ is split semisimple where $F=F_0(v)$ (see \cite[9.3.5]{ourbuch}). We 
now take the ring
\[ R=\{f/g \mid f,g \in R_0[v,v^{-1}], g\not\in \fp[v,v^{-1}]\}\supseteq 
A_p\supseteq \cO.\]
Then $R$ is a discrete valuation ring with field of fractions $F$ and
residue field $k$. The hypotheses in (\ref{remcent})(0) are satisfied
and we have a decomposition matrix $D_R$.
\end{abschnitt}

\begin{cor} \label{cc} Assume that there exists a discrete valuation $R$
as in (\ref{remcent}) which is $\bc$-adapted for every $\bc\in B_0$. Then 
the set  
\[ \Big\{\sum_{E\in \Irr(W)} [E:[\Gamma]]\, E \; \Big|\; 
\mbox{ $\Gamma$ any left cell of $W$}\Big\} \subseteq {\Q}[\Irr(W)]\]
is linearly independent. 
\end{cor} 

\begin{proof} This is just another way of saying that the columns of
the decomposition matrix $D_R$ are linearly independent. \end{proof}

\begin{rem} \label{cc1} If $W$ is of classical type, then $R=A_2$ 
satisfies the above hypothesis; see Section~6. For groups of 
exceptional type, a ring $R$ satisfying the hypothesis of Corollary~\ref{cc}
no longer exists. However, one can just check by an explicit computation 
(using the results in Section~7) that the above statement holds. Thus, 
the statement in Corollary~\ref{cc} is seen to hold for any $W,L$ 
(assuming (P1)--(P15)).
\end{rem}

\section{Constructible representations and families} \label{a-func}

We preserve the set-up of the previous sections. Now we turn to
constructible representations and families. The definition relies on
Lusztig's notion of {\em truncated induction}.

First recall the definition of the invariants $\ba_E$ for any $E\in \Irr(W)$;
see (\ref{simpH}). Now, given any ${\C}[W]$-module $E'$, we can write
uniquely $E'\cong E_0'\oplus E_1'\oplus E_2'\oplus \cdots$ where, for any 
integer $i$, we define 
\[ E_a':=\bigoplus_E [E:E']\, E \quad \mbox{(sum over all $E\in \Irr(W)$ 
such that $\ba_E=i$)};\]
here, $[E:E']$ denotes the multiplicity of $E$ as a simple component of
$E'$. Now let $I\subset S$ and consider the parabolic subgroup $W_I 
\subset W$. Let $E'$ be any ${\C}[W_I]$-module and denote by $\Ind_{I}^S(E')$ 
the corresponding ${\C}[W]$-module obtained by induction. Then, as above, we 
can write uniquely 
\[ \Ind_{I}^S(E')\cong E_0 \oplus E_1 \oplus E_2 \oplus \cdots,\quad 
\mbox{where} \quad E=\Ind_I^S(E').\]
It is shown in \cite[20.15]{Lusztig03} that if $E'=E_a' \neq \{0\}$ for 
some $a$, then we have 
\[ \Ind_{I}^S(E')\cong E_a \oplus E_{a+1} \oplus E_{a+2} \oplus\cdots,
\qquad \mbox{where $E_a\neq \{0\}$}.\]
In this case, we set $\Jind_I^S(E'):=E_a$. Thus, $\Jind_{I}^S(E')$ is a
${\C}[W]$-module such that 
\[\Ind_{I}^S(E') \cong\Jind_{I}^S(E')\oplus\mbox{ higher terms},\]
where ``higher terms'' means a direct sum of simple $W$-modules with 
$\ba$-invariant strictly bigger than $a$. This is precisely Lusztig's 
{\em truncated induction} in \cite[20.15]{Lusztig03}.

\begin{abschnitt} {\bf Constructible representations and families.} 
\label{fam} Following Lusztig \cite[22.1]{Lusztig03}, we define the set 
$\Con(W)$ of constructible representations of $W$ (with respect to the 
fixed weight function $L$). This is done by induction on $|W|$. If 
$W=\{1\}$, then $\Con(W)$ consists of the unit representation. If $W\neq 
\{1\}$, then $\Con(W)$ consists of the ${\C}[W]$-modules  of the form 
$\Jind_{I}^S(E')$ or $\Jind_{I}^S(E')\otimes \mbox{sgn}$, for various 
subsets $I \subsetneqq S$ and various $E'\in \Con(W_I)$. 

We define a corresponding ``decomposition matrix'' $\bD$ as follows. The 
rows are labelled by $\Irr(W)$ and the columns are labelled by $\Con(W)$;
the coefficients in a fixed column of $\bD$ give the expansion of the
corresponding constructible representation in terms of the irreducible
ones. 

The matrix $\bD$ has been computed explicitly (for all $W,L$) by Lusztig 
\cite{Lusztig82b}, \cite[Chap.~22]{Lusztig03} and Alvis--Lusztig 
\cite{AlLu} (type $H_4$).

Following Lusztig \cite[23.1]{Lusztig03}, we define a graph $\cG_W$ as 
follows. The vertices of $\cG_W$ are labelled by $\Irr(W)$. Given $E\neq E'$
in $\Irr(W)$, the corresponding vertices in $\cG_W$ are joined if $E,E'$
both appear as simple components of some constructible representation
of $W$. We say that $E,E' \in \Irr(W)$ belong to the same {\em family} if
the corresponding vertices are in the same connected component of $\cG_W$.

The partition of $\Irr(W)$ and $\Con(W)$ according to families gives rise
to a block diagonal shape of $\bD$, with one block on the diagonal for each
family.
\end{abschnitt}

\begin{rem} \label{con1} One may also consider the following set 
$\Con'(W)$ of representations of $W$ (with respect to the fixed weight 
function $L$). The definition of $\Con'(W)$ appears, in a slightly 
refined form, in Malle--Rouquier \cite[\S 2]{MaRo}.

Again, we proceed by induction on $|W|$. If $W=\{1\}$, then $\Con'(W)$ 
consists of the unit representation. If $W\neq \{1\}$, then $\Con'(W)$ 
consists of all ${\C}[W]$-modules of the form $E_i$, where $E=\Ind_I^S(E')$,
for various subsets $I \subseteq S$ such that $|I|=|S|-1$, various 
$E'\in \Con'(W_I)$ and various $i\in \Z$. It is clear that 
\[ \Con(W) \subseteq \Con'(W).\]
Note that we certainly do not have $\Con(W)=\Con'(W)$ in general. For 
example, for $W$ of type $E_7$, we have $2\cdot (512_a\oplus 512_a') \in 
\Con'(W)\setminus \Con(W)$ (see the discussion in (\ref{celle7}) below). 
In Remark~\ref{con1a}, we shall see that $\Con'(W) \subseteq {\N}[\Con(W)]$ 
\end{rem} 

\begin{thm}[Lusztig, Rouquier] \label{thmfam} Let $\cO\subseteq {\C}(v)$ 
be the subring defined in Example~\ref{exp1}(b). Assume that (P1)--(P15) 
hold for $W,L$. Let $E,E' \in\Irr(W)$.  Then the following conditions 
are equivalent.
\begin{itemize}
\item[(a)] $E,E'$ belong to the same family, in the sense of (\ref{fam}).
\item[(b)] $E,E'$ belong to the same two-sided cell, that is, there
exists some $\bc \in B_0$ such that $t_{\bc}E_{\spadesuit}\neq 0$ and 
$t_{\bc}E_{\spadesuit}'\neq 0$; see \cite[20.2]{Lusztig03}.
\item[(c)] $E,E'$ belong to the same block of $H_{\cO}$, that is, there 
exists a primitive idempotent $e$ in the center of $H_{\cO}$ such that 
$eE_v\neq 0$ and $eE_v'\neq 0$.
\end{itemize}
Furthermore, if $E,E'$ belong to the same family, then $\ba_E=\ba_{E'}$.
\end{thm}

\begin{proof} For the equivalence of (a) and (b), see the argument in the
proof of \cite[Prop.~23.3]{Lusztig03}. The equivalence  of (b) and (c) is an
immediate consequence of Theorem~\ref{prop32}(b). For the statement 
concerning $\ba_E$ and $\ba_{E'}$; see \cite[Prop.20.6(c)]{Lusztig03}.
\end{proof}

The equivalence of (a) and (b) was first proved by Barbasch and Vogan in
the equal parameter case; see \cite[Theorem~5.25]{LuBook}. The equivalence 
of (b) and (c) has been established by Rouquier \cite{rou} in the equal
parameter case. The important point about that equivalence is that the
statement in (c) is independent of the Kazhdan--Lusztig basis $\{c_w\}$ of 
$H$. As such, it applies to a wider class of algebras; see Brou\'e--Kim 
\cite{BK} and Malle--Rouquier \cite{MaRo}.

We now establish some results which will be helpful for the computation of
the left cell representations. The discussion mainly follows 
\cite{Lusztig03}.  Note that the {\em statements} in
Lemmas~\ref{corfam} and~\ref{lem12} (and the corollaries following them)
do not involve the properties (P1)--(P15): these are only implicit in
the proofs. It would be very interesting to prove any of those 
statements by elementary methods.

\begin{lem} \label{corfam} Let $\Gamma$ be a left cell of $W$. Then
all simple components of $[\Gamma]$ belong to a fixed family of $\Irr(W)$.
In particular, all  simple components of $[\Gamma]$ have the same
$\ba$-invariant.
\end{lem}

\begin{proof} Let $\cD\cap \Gamma=\{d\}$. If $E,E'$ are simple components 
of $[\Gamma]$, then $t_d$ acts non-trivially on $E_{\spadesuit}$ and 
on $E_{\spadesuit}'$; see \cite[Prop.~21.4]{Lusztig03}. Consequently,
we also have $t_{\bc}E_{\spadesuit}\neq 0$ and $t_{\bc}E_{\spadesuit}'
\neq 0$, where $\bc$ is the two-sided cell containing~$d$. Thus, 
$E,E'$ belong to the same two-sided cell. Now the assertion follows 
from Theorem~\ref{thmfam}.
\end{proof}

\begin{lem} \label{lem10} Let $\Gamma$ be a left cell of $W$. Then
$\Gamma w_0$ also is a left cell of $W$, where $w_0$ is the unique element 
of maximal length of $W$. We have $[\Gamma w_0]\cong [\Gamma]\otimes 
\operatorname{sgn}$ where $\operatorname{sgn}$ is the sign 
representation of $W$.
\end{lem}

\begin{proof} See \cite[Cor.~11.7]{Lusztig03} and 
\cite[Prop.~21.5]{Lusztig03}. Here, (P1)--(P15) are not needed.
\end{proof}

The following result is an extremely strong condition on
the expansion of a left cell representation in terms of the
irreducible ones.

\begin{lem} \label{lem12} Let $\Gamma$ be a left cell of $W$. Then 
we have 
\[ 1=\sum_{E\in \Irr(W)} \frac{1}{f_E}\, [E:[\Gamma]].\]
Here, $[E:[\Gamma]]$ denotes the multiplicity of $E$ 
as a simple component of $[\Gamma]$.
\end{lem}

\begin{proof} Let $\Gamma\cap\cD=\{d\}$ and $n_d=\pm 1$. By 
\cite[19.2 and 20.1(b)]{Lusztig03}, we have
\[ n_d=\tau(t_d)=\sum_{E\in \Irr(W)} \frac{1}{f_E}\, 
\mbox{tr}(t_d,E_\spadesuit).\]
On the other hand, by \cite[Prop.~21.4]{Lusztig03}, we have $n_d
\mbox{tr}(t_d,E_\spadesuit)=[E:[\Gamma]]$ for every $E\in \Irr(W)$.
This yields the desired assertion.
\end{proof}

\begin{cor} \label{cor12} Let $\Gamma$ be a left cell of $W$. Assume
that there exists some $E\in \Irr(W)$ such that $[E:[\Gamma]]\neq 0$
and $f_E=1$. Then $[\Gamma]\cong E\in \Irr(W)$ and $[\Gamma]$ is 
constructible.
\end{cor}

\begin{proof} Let $\Irr(W)=\{E_1,\ldots,E_r\}$ where $E_1=E$. We write
$[\Gamma]\cong \bigoplus_{i=1}^r n_i\, E_i$ where $n_i\in {\Z}_{\geq 0}$. 
Then, by Lemma~\ref{lem12}, we have 
\[ 1=\sum_{i=1}^r \frac{1}{f_{E_i}}\, n_i=
n_1+\sum_{i=2}^r \frac{1}{f_{E_i}}\, n_i.\]
Now note that $n_1=[E:[\Gamma]]$ is a positive integer and that $f_{E_i}>0$
for all $i$. Hence we must have $n_1=1$ and $n_i=0$ for all $i>1$.

To prove that $E\cong [\Gamma]$ is constructible, we argue as follows.
By \cite[Prop.~22.3]{Lusztig03}, there exists some constructible
representation $E'$ of $W$ which contains a simple component isomorphic to
$E$. By \cite[Lemma~22.2]{Lusztig03}, there exists a left cell $\Gamma'$
such $E'\cong [\Gamma']$. Thus, $E$ is a simple component of $[\Gamma']$.
By the previous argument, we have $E\cong [\Gamma']$ and so $E\cong E'$, 
as desired.
\end{proof}

\begin{cor} \label{cor12a} Assume that $f_E=1$ for all $E\in \Irr(W)$.
Then $[\Gamma]$ is irreducible and constructible for every left cell 
$\Gamma$. 
\end{cor}

\begin{proof} Immediate by Corollary~\ref{cor12}. 
\end{proof}

\begin{exmp} \label{expA} Assume that $W$ is of type $A_{m-1}$ ($m\geq 2$) 
with generators and relations given by the following diagram:
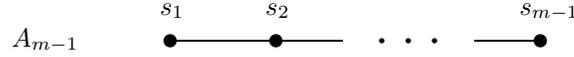
\begin{center}
\begin{picture}(300,30)
\put( 30, 8){$A_{m-1}$}
\put( 90 ,10){\circle*{5}}
\put(130 ,10){\circle*{5}}
\put(230 ,10){\circle*{5}}
\put( 90, 10){\line(1,0){40}}
\put(130, 10){\line(1,0){25}}
\put(170, 10){\circle*{2}}
\put(180, 10){\circle*{2}}
\put(190, 10){\circle*{2}}
\put(205, 10){\line(1,0){25}}
\put( 86 ,20){$s_1$}
\put(126, 20){$s_2$}
\put(222, 20){$s_{m{-}1}$}
\end{picture}
\end{center}
All generators are conjugate in $W$. So every weight function  is of 
the form $L=al$ where $a>0$. We have $f_E=1$ for all $E \in \Irr(W)$; see 
\cite[22.4]{Lusztig03}. Now Corollary~\ref{cor12a} shows that $[\Gamma]$ 
is irreducible and constructible for every left cell $\Gamma$. On the other 
hand, by \cite[22.5]{Lusztig03}, the irreducible representations are
precisely the constructible ones. The fact that $[\Gamma]$ is irreducible
for every left cell $\Gamma$ has been already established by
Kazhdan--Lusztig \cite[\S 5]{KL}; for a more  detailed and completely 
elementary proof, see Ariki \cite{Ar2}.
\end{exmp}

\begin{exmp} \label{expB1} Let $W$ be of type $B_m$ ($m \geq 2$), with 
generators and relations given by the following diagram:
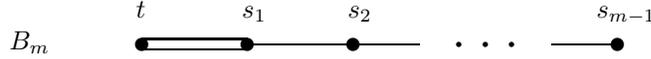
\begin{center}
\begin{picture}(300,30)
\put(  0, 8){$B_m$}
\put( 50 ,10){\circle*{5}}
\put( 50,  8){\line(1,0){40}}
\put( 50, 12){\line(1,0){40}}
\put( 48, 20){$t$}
\put( 90 ,10){\circle*{5}}
\put(130 ,10){\circle*{5}}
\put(230 ,10){\circle*{5}}
\put( 90, 10){\line(1,0){40}}
\put(130, 10){\line(1,0){25}}
\put(170, 10){\circle*{2}}
\put(180, 10){\circle*{2}}
\put(190, 10){\circle*{2}}
\put(205, 10){\line(1,0){25}}
\put( 88 ,20){$s_1$}
\put(128, 20){$s_2$}
\put(222, 20){$s_{m{-}1}$}
\end{picture}
\end{center}
Let $L$ be a weight function such that
\[ L(t)=b>0 \qquad \mbox{and}\qquad L(s_1)=L(s_2)=\cdots =L(s_{m-1})=a>0.\]
Assume that (P1)--(P15) hold for $W,L$. As in \cite[22.6]{Lusztig03}, we 
write $b=ar+b'$ where $r,b'\in \Z$ such that $r\geq 0$ and $0\leq b'<a$.

Assume first that $b'>0$. Then we have $f_E=1$ for all $E\in \Irr(W)$;
see \cite[Prop.~22.14(a)]{Lusztig03}. So Corollary~\ref{cor12a} shows
that $[\Gamma]\in \Irr(W)$ for all left cells $\Gamma$. On the other
hand, the constructible representations are precisely the irreducible
ones by \cite[Prop.~22.25]{Lusztig03}. In particular, Conjecture~\ref{conjC}
holds for $W,L$ if $b'>0$.
The case where $b'=0$ will be considered in Section~6.
\end{exmp}

\begin{rem} \label{rem3} Let $W$ be of type $B_m$ with parameters
given as above. Assume that $b/a$ is ``large''. Then the left cells are 
explicitly determined by Bonnaf\'e--Iancu \cite{BI} (without using 
the assumption that (P1)--(P15) hold); the corresponding representations 
are irreducible and constructible. There is some hope that similar arguments
might be found to deal with arbitrary values of $a$ and $b$, as long as 
$a$ does not divide~$b$.
\end{rem}

Finally, we present some further conditions on the coefficients in the 
expansion $[\Gamma]\cong \bigoplus_E [E:\Gamma]]\, E$ for a left cell 
$\Gamma$. This will also be needed in the discussion of type $E_8$ in
(\ref{celle8}). The desired conditions follow from results in \cite{GeRo} 
concerning the center of $H$; we shall follow the exposition in 
\cite{ourbuch}. 

\begin{abschnitt} {\bf Central characters.} \label{center} Let $Z(H)$
be the center of $H$. An element $z\in Z(H)$ acts by a scalar in
every simple $H_K$-modules $E_v$, where $E \in \Irr(W)$; we denote that
scalar by $\omega_E(z)$. For technical simplicity, let us now
assume that $W$ is a Weyl group; then $\Q$ is a splitting field for
$W$ and ${\Q}(v)$ is a splitting field for $H$; see 
\cite[Thoerem~9.3.5]{ourbuch}. Consequently, since $A$ is integrally
closed in ${\Q}(v)$, we have 
\begin{equation*}
\omega_E(z)\in A\qquad \mbox{for any $E\in \Irr(W)$ and $z\in Z(H)$}.
\tag{a}
\end{equation*}
Now let $R \subseteq {\C}(v)$ be a noetherian subring such that 
$\cO\subseteq R$, where $\cO$ is the ring in Example~\ref{exp1}(b). 
Let $e\in H_R$ be an indempotent and consider the corresponding 
projective $H_R$-module $P:=H_Re$. Extending scalars from $R$ to
$K$, we obtain an $H_K$-module $P_K$. We denote 
\[ [E:P]=[E_v:P_K]:=\mbox{multiplicity of $E_v$ in $P_K$, for any 
$E\in \Irr(W)$}.\] 
Then the argument in the proof of
\cite[Theorem~7.5.3]{ourbuch} shows that 
\begin{equation*}
\sum_{E\in \Irr(W,\bc)} \frac{1}{c_E} \,[E:P]\,\omega_E(z_C) 
\in R \qquad \mbox{for any $z\in Z(H)$}.\tag{b}
\end{equation*}
This yields rather restrictive conditions on the coefficients $[E:P]$.

By Corollary~\ref{lem34}, we can apply this, in particular, to the module 
$P=[\Gamma]_R$ for a left cell $\Gamma$ of $W$. Note that, for $z=1$, 
we have $\omega_E(1)=1$ for all $E$ and so the above condition bears 
some resemblance to that in Lemma~\ref{lem12}.

In order to be able to use the formula (b), we shall need to compute 
$\omega_E(z)$ for some elements $z\in Z(H)$.

Now, by \cite[Cor.~8.2.5]{ourbuch}, there is a distinguished $A$-basis 
of $Z(H)$, denoted by $\{z_C\mid C \in \Cl(W)\}$ where $\Cl(W)$ is
the set of conjugacy classes of $W$. The scalars $\omega_E(z_C)$ are 
determined by the following identity:
\begin{equation*}
\sum_{C \in \Cl(W)} v_w\,\mbox{tr}(T_{w_C},E_v)\,\omega_{E'}(z_C)=
\left\{\begin{array}{cl} c_E & \quad \mbox{if $E\cong E'$},\\ 0 & \quad
\mbox{otherwise};\end{array}\right.\tag{c}
\end{equation*}
see \cite[Exc.~9.5]{ourbuch}. Here, $w_C$ is an element of minimal 
length in $C\in \Cl(W)$. Note that, by \cite[Cor.~8.2.6]{ourbuch}, the
value $\mbox{tr}(T_{w_C},E_v)$ does not depend on the choice of $w_C$.
The ``character tables'' 
\[ \Big(\mbox{tr}(T_{w_C},E_v)\Big)_{E\in \Irr(W), C\in \Cl(W)}\]
are explicitly known for all $W,L$; see \cite{ourbuch}. As explained in 
the proof of \cite[Prop.~11.5.13]{ourbuch}, the identities (c) can be
used to compute the scalars $\omega_E(z_C)$. 
\end{abschnitt}

\section{Induction and restriction of left cells} \label{indres}

Let $I\subseteq S$ and consider the parabolic subgroup $W_I$. The restriction
of $L$ to $W_I$ is a weight function on $W_I$. Thus, we have a partition 
of $W_I$ into left cells with respect to $L|_{W_I}$. We shall now consider 
the compatibility of the corresponding left cell representations with  
respect to induction and restriction. All the results that we are going
to present here were already known in the case of equal parameters. However,
new arguments are needed in the proofs for some of these results in the 
unequal parameter case (most notably Lemmas~\ref{lem42}, \ref{propJ} and 
\ref{bv2}).

\begin{lem}[See Barbasch--Vogan \protect{\cite[Prop.~3.11]{BV}} in the 
equal parameter case] \label{lem41} Let $\Gamma$ be a left cell of $W$ 
and $I \subset S$. Then there exist (pairwise different) left cells 
$\Gamma_1',\ldots,\Gamma_r'$ of $W_I$ and positive integers $n_1,\ldots,
n_r$ such that 
\[ {\Res}_{I}^S([\Gamma])\cong n_1 [\Gamma_1']\oplus \cdots 
\oplus n_r[\Gamma_r'].\]
\end{lem}

\begin{proof} We follow the proof given by Roichman
\cite[Theorem~5.2]{Roich} in the equal parameter case. 

Recall the definition of the relation $\leq_{\cL}$ on $W$: this is the 
transitive closure of the relation ``$y \leftarrow_{\cL} w$ if $h_{s,w,y}
\neq 0$ for some $s\in S$''; see \cite[8.1]{Lusztig03}. 

We define a relation $\leq_{\cL,I}$ on $W$ as the transitive
closure of the relation ``$y \leftarrow_{\cL,I} w$ if $h_{s,w,y}\neq 0$ 
for some $s\in I$''. Let $\sim_{\cL,I}$ be the corresponding equivalence
relation on $W$. The restriction of $\leq_{\cL,I}$ to $W_I$ gives 
precisely the left cells of $W_I$ with respect to $L|_{W_I}$.

Let $Y_I$ be the set of all $w\in W$ such that $w$ has minimal length 
in the right coset $W_Iw$. Let $u,u' \in W_I$ and $x,x'\in Y_I$.  
Now, if $ux\leftarrow_{\cL,I} u'x'$, then there exists some $s\in I$ such 
that $h_{s,u'x',ux}\neq 0$. By the formula in \cite[Prop.~6.3]{Lusztig03},
we have either $ux=su'x'>u'x'$ or $ux<u'x'$. Thus, we have the implication
\[ ux \leq_{\cL,I} u'x'\quad \mbox{and}\quad l(ux)>l(u'x') \quad 
\Rightarrow \quad x=x' .\]
Consequently, we have the implication
\[ ux \sim_{\cL,I} u'x' \quad \Rightarrow \quad x=x'.\]
This shows that there exist (pairwise different) left cells $\Gamma_1',
\ldots,\Gamma_r'$ of $W_I$ and subsets $R_1,\ldots,R_r$ of $Y_I$ such that
\[\Gamma=\bigcup_{i=1}^r \bigcup_{y\in R_i} \Gamma_i'y.\]
Now consider the restriction of the $H$-module $[\Gamma]_A$ to $H_I$.
The formula in \cite[Lemma~9.10(e)]{Lusztig03} shows that, for
fixed $i \in \{1,\ldots,r\}$ and $y\in R_i$, we have 
\[ c_s^\dagger \cdot c_{uy}^\dagger =\sum_{u'\in \Gamma_i'} h_{s,u,u'}
c_{u'y}^\dagger \qquad \mbox{for any $s\in I$ and $u\in \Gamma_i'$}.\]
Thus, as an $H_I$-module, we have $\Res_{I}^S([\Gamma]_A)\cong n_1\, 
[\Gamma_1']_A\oplus \cdots \oplus n_r\,[\Gamma_r']_A$, where $n_i=|R_i|$
for all $i$.  Upon setting $v\mapsto 1$, we obtain the required assertion 
concerning $W_I$-modules.
\end{proof}

\begin{lem}[See Barbasch--Vogan \protect{\cite[Prop.~3.15]{BV}} in the 
equal parameter case] \label{lem42} Let $I\subset S$ and $\Gamma'$ be a 
left cell of $W_I$. Then there exist (pairwise different) left cells 
$\Gamma_1,\ldots,\Gamma_r$ of $W$ such that 
\[ {\Ind}_{I}^S([\Gamma'])\cong [\Gamma_1]\oplus \cdots \oplus [\Gamma_r].\]
\end{lem}

\begin{proof} Let $X_I$ be the set of all $x\in W$ such that $x$ has 
minimal length in the left coset $xW_I$. Then, by \cite{myind}, we have
$X_I\Gamma'=\Gamma_1\amalg \cdots \amalg \Gamma_r$ where $\Gamma_i$ are 
left cells of $W$.  So we have an $H_A$-module
\[ [X_I\Gamma']_A:=\sum_{x \in X_I}\sum_{u\in \Gamma'} Ac_{xu}^\dagger
\cong [\Gamma_1]_A\oplus \cdots \oplus [\Gamma_r]_A.\]
Here, $c_w^\dagger$ ($w\in W$) acts by the rule $c_w^\dagger \cdot 
c_y^\dagger=\sum_{z\in X_I\Gamma'} h_{x,y,z}c_z^\dagger$ for any $y\in 
X_I\Gamma'$. 
To identify this module with an induced module, we set
\begin{align*}
\cI &:=\langle T_xc_u^\dagger \mid x\in X_I,u\in W_I,u \leq_{\cL,I} u' 
\mbox{ for some $u' \in \Gamma'$}\rangle_A,\\
\hat{\cI} &:=\langle T_xc_u^\dagger \mid x\in X_I, u \in W_I\setminus
\{\Gamma'\}, u \leq_{\cL,I} u' \mbox{ for some $u' \in \Gamma'$}\rangle_A.
\end{align*}
Then $\hat{\cI} \subseteq \cI$ are left ideals in $H_A$; see
\cite[Lemma~2.2]{myind}. Thus, $\cI/\hat{\cI}$ is an $H_A$-module; it
is free as an $A$-module with a basis given the residue classes of the
elements $T_xc_u^\dagger$ ($x\in X_I$, $u\in \Gamma'$). By the definition 
of $\Ind_I^S$ (see \cite[\S 9.1]{ourbuch}), we have  an isomorphism of
$H_A$-modules 
\[ \Ind_I^S([\Gamma']) \stackrel{\sim}{\longrightarrow} \cI/\hat{\cI},
\qquad T_x\otimes c_u^\dagger \mapsto T_xc_u^\dagger +\hat{\cI}.\]
On the other hand, by \cite[Prop.~3.3]{myind}, we also have 
\begin{align*}
\cI &=\langle c_{xu}^\dagger \mid x\in X_I,u\in W_I,u \leq_{\cL,I} u' 
\mbox{ for some $u' \in \Gamma'$}\rangle_A,\\
\hat{\cI} &=\langle c_{xu}^\dagger \mid x\in X_I, u \in W_I\setminus
\{\Gamma'\}, u \leq_{\cL,I} u' \mbox{ for some $u' \in \Gamma'$}\rangle_A.
\end{align*}
Thus, $\cI/\hat{\cI}$ also has an $A$-basis given by the residue 
classes of the elements $c_{xu}^\dagger$ ($x\in X_I$, $u\in \Gamma'$). Hence 
there is an isomorphism of $H_A$-modules
\[ [X_I\Gamma']_A \stackrel{\sim}{\longrightarrow} \cI/\hat{\cI},
\qquad c_{xu}^\dagger \mapsto c_{xu}^\dagger +\hat{\cI}.\]
Upon setting $v\mapsto 1$, we obtain $\Ind_I^S([\Gamma']) \cong [\Gamma_1]
\oplus \cdots \oplus [\Gamma_r]$ as desired.
\end{proof}

\begin{rem} \label{donot} The proofs of Lemma~\ref{lem41} and
Lemma~\ref{lem42} do not require the assumption that (P1)--(P15) hold.
\end{rem}

\begin{rem} \label{con1a} We can now clarify the relation between 
$\Con(W)$ and the set $\Con'(W)$ defined in Remark~\ref{con1}. (See 
Malle--Rouquier \cite[Prop.~2.5]{MaRo} for the case of equal parameters.)
Assume that Conjecture~\ref{conjC} holds for $W,L$ and all parabolic 
subgroups of $W$. Then we have:
\begin{itemize}
\item[(a)] $\Con'(W) \subseteq {\N}[\Con(W)]$ and
\item[(b)] the regular representation of $W$ is a sum of constructible
representations.
\end{itemize}
Indeed, using the notation in Remark~\ref{con1}, let $E_i\in \Con'(W)$ 
where $E=\Ind_I^S(E')$ for some $E' \in \Con'(W_I)$ ($I\subsetneqq S$) 
and some integer $i\geq 0$. By induction, $E'$ is a sum of constructible
representations. Hence, by Conjecture~\ref{conjC}, $E'$ is a sum of left 
cell representations. Using Lemma~\ref{corfam} and Lemma~\ref{lem42}, we 
conclude that $E_i$ is a sum of left cell representations and, hence 
(via Conjecture~\ref{conjC}), a sum of constructible representations.
This proves (a). Now (b) is an easy consequence: we just have to note 
that the regular representation of $W$ (which can be obtained by inducing
the regular representation of any proper parabolic subgroup) is a sum
of representations in $\Con'(W)$.
\end{rem}

\begin{lem} \label{propJ} Let $I\subset S$ and $\Gamma'$ be a left 
cell of $W_I$. Then we have 
\[ \Jind_{I}^S([\Gamma'])\cong [\Gamma],\]
where $\Gamma$ is the left cell of $W$ such that $\Gamma'\subseteq \Gamma$.
\end{lem}

\begin{proof} See the argument in the proof of Case~1 in 
\cite[Lemma~22.2]{Lusztig03}. Note that this argument is much simpler
than one given in \cite[5.28.1]{LuBook} in the case of equal parameters,
which uses results from Barbasch--Vogan \cite{BV}.
\end{proof}

\begin{lem}[See Lusztig \protect{\cite[\S 3]{Lusztig86}} in the equal
parameter case] \label{bv2} Let $\cF\subseteq \Irr(W)$ be a family and 
assume that there exists a proper subset $I \subset S$ and a family 
$\cF'\subset \Irr(W_I)$ such that $\Jind_{I}^S$ defines a bijection 
$\cF'\stackrel{\sim}{\rightarrow} \cF$.  

Let $\Gamma$ be a left cell of $W$ and assume that all simple components 
of $[\Gamma]$ belong to $\cF$. Then there exists a left cell $\Gamma'$ of 
$W_I$ such that $[\Gamma] \cong\Jind_I^S([\Gamma'])$. In particular, if
$[\Gamma']$ is constructible, then so is $[\Gamma]$.
\end{lem}

\begin{proof} Let $\cF=\{E_1,\ldots,E_f\}$. By assumption, we have
$\cF'=\{E_1',\ldots,E_f'\}$ where $E_i\cong \Jind_I^S(E_i')$ for all~$i$.
Using \cite[Prop.~20.13]{Lusztig03}, we conclude that  
\[ \Res_I^S(E_i)\cong E_i'\; \oplus\mbox{ combination of $E'\in \Irr(W_I)$ 
with $\ba_{E'}< a$}\]
where $a=\ba_{E_1}=\cdots =\ba_{E_f}$. Now let $\Gamma$ be a left cell 
belonging to $\cF$; we have $[\Gamma]\cong  \bigoplus_i[E_i:[\Gamma]]\, 
E_i$. Upon restriction to $W_I$, we obtain
\[\Res_I^S([\Gamma])\cong \bigoplus_{i=1}^f [E_i:[\Gamma]]\, E_i'\;
 \oplus \mbox{ combination of $E'$ with $\ba_{E'}<a$}.\]
On the other hand, by Lemma~\ref{lem41}, we know that $\Res_{I}^S([\Gamma])$
is a sum of left cell representations of $W_I$. Since all simple 
components of a left cell representation belong to a fixed family, we 
conclude that 
\[\bigoplus_{i=1}^f [E_i:[\Gamma]]\, E_i'\; \cong \mbox{ sum of left cell 
representations of $W_I$}.\]
Applying $\Jind_I^S$ and using Lemma~\ref{propJ}, we deduce that
\[[\Gamma]\cong \bigoplus_{i=1}^f [E_i:[\Gamma]]\, E_i \; \cong \mbox{ sum 
of left cell representations of $W$},\]
where each representation occuring in the sum is obtained by 
$\Jind$-induction of a left cell representation of $W_I$. Thus, to complete 
the proof, it is sufficient to show that only one left cell representation
occurs in the above sum. But this follows from the following general
statement:

\medskip
{\em Let $\Gamma,\Gamma_1$ be left cells and assume that $[\Gamma_1]\cong 
[\Gamma]\oplus E'$ for some ${\C}[W]$-module $E'$. Then we automatically 
have $[\Gamma]\cong [\Gamma_1]$}.
\medskip

Indeed, the assumption implies that $[E:[\Gamma]]\leq [E:\Gamma_1]]$
for any $E \in \Irr(W)$. Hence Lemma~\ref{lem12} yields 
\[ 1=\sum_{E \in \Irr(W)} \frac{1}{f_E}\, [E:[\Gamma]]
\leq \sum_{E \in \Irr(W)} \frac{1}{f_E}\, [E:[\Gamma_1]]=1.\]
Consequently, all the inequalities must be equalities and so
$[\Gamma]\cong [\Gamma_1]$.
\end{proof}

Note that the proof of Lemma~\ref{bv2} which is given in 
\cite[\S 3]{Lusztig86}, uses the fact that every left cell representation 
has a unique ``special'' simple component with multiplicity~$1$; see 
the discussion in \cite[5.25]{LuBook}.  This argument is no longer 
available in the case of unequal parameters: there are families which do
not contain any reasonably defined ``special'' simple module (see, 
for example, \cite[Remark~4.11]{mykl}).

\begin{abschnitt} {\bf Cuspidal families.} \label{cusp} Let
$\cF$ be a family of $\Irr(W)$. Then $\cF \otimes \mbox{sgn}=
\{ E\otimes \mbox{sgn} \mid E \in \cF\}$ also is a family. (This is
clear by the definitions in (\ref{fam}).) As in \cite[8.1]{LuBook},
we say that $\cF$ is {\em cuspidal} if neither $\cF$ nor $\cF 
\otimes \mbox{sgn}$ satisfies the hypothesis of Lemma~\ref{bv2}
for any proper subset $I \subset S$ and any family $\cF'$ of $\Irr(W_I)$.

Suppose now we want to show that $[\Gamma]$ is constructible for
every left cell $\Gamma$ of $W$. Proceeding by induction on $|S|$,
we may assume that $[\Gamma']$ is constructible for every left cell
$\Gamma'$ of a proper parabolic subgroup of $W$. Then, as explained
in \cite{Lusztig86}, it will be enough to consider only those left 
cells $\Gamma$ for which $[\Gamma]$ belongs to a cuspidal family of 
$\Irr(W)$. 

Indeed, let $\Gamma$ be a left cell such that all simple components of 
$[\Gamma]$ belong to a family $\cF$ which is not cuspidal. Then either 
$\cF$ itself or $\cF \otimes \mbox{sgn}$ satisfies the hypothesis of 
Lemma~\ref{bv2} for some proper subset $I \subset S$ and some family 
$\cF'$ of $\Irr(W_I)$. Assume first that $\cF$ itself satisfies that
hypothesis. Then there exists a left cell $\Gamma'$ of $W_I$ such that
$[\Gamma]\cong \Jind_I^S([\Gamma'])$. By induction, $[\Gamma']$ is 
constructible, hence $[\Gamma]$ is construcible.  Now assume that 
$\cF \otimes\mbox{sgn}$ satisfies that hypothesis. Then consider the 
left cell $\Gamma w_0$. By Lemma~\ref{lem10}, all simple components 
of $[\Gamma w_0]\cong [\Gamma]\otimes \mbox{sgn}$ belong to $\cF \otimes 
\mbox{sgn}$.  Hence, as before, we may conclude that $[\Gamma w_0]$ 
is constructible. Consequently, $[\Gamma]$ also is constructible. 
\end{abschnitt}

\section{Types $B_m$ and $D_m$} \label{bdcells}
We are now going to prove Conjecture~\ref{conjC} for groups of type $B_m$ 
and $D_m$, under the hypothesis that (P1)--(P15) hold for $W,L$. A special
case in type $B_m$ has been already considered in Example~\ref{expB1}.
(For type $A_{m-1}$, see Example~\ref{expA}.) As far as type $D_m$ and $B_m$
with equal parameters are concerned, our proof is different from the one 
given by Lusztig \cite{Lusztig86}.

\begin{abschnitt} {\bf A decomposition matrix in type $B_m$.} \label{typeB} 
Let $W$ be of type $B_m$ ($m \geq 2$), with generators and relations 
given by the following diagram:
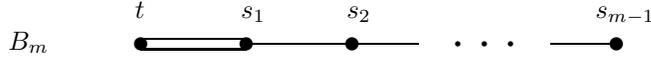
\begin{center}
\begin{picture}(300,30)
\put(  0, 8){$B_m$}
\put( 50 ,10){\circle*{5}}
\put( 50,  8){\line(1,0){40}}
\put( 50, 12){\line(1,0){40}}
\put( 48, 20){$t$}
\put( 90 ,10){\circle*{5}}
\put(130 ,10){\circle*{5}}
\put(230 ,10){\circle*{5}}
\put( 90, 10){\line(1,0){40}}
\put(130, 10){\line(1,0){25}}
\put(170, 10){\circle*{2}}
\put(180, 10){\circle*{2}}
\put(190, 10){\circle*{2}}
\put(205, 10){\line(1,0){25}}
\put( 88 ,20){$s_1$}
\put(128, 20){$s_2$}
\put(222, 20){$s_{m{-}1}$}
\end{picture}
\end{center}
Let $L$ be a weight function such that 
\[L(t)=r\geq 0\qquad\mbox{and}\qquad L(s_1)=L(s_2)=\cdots =L(s_{m-1})=1.\]
Here, we explicitly allow the case where $r=0$, which is related to 
type $D_m$; see (\ref{typedn}). The definition of left cells and
constructible representations still makes sense in this case; see 
\cite{Lusztig03}. Let us assume that (P1)--(P15) hold for $W,L$. It is
known that this is the case if $r=1$ (equal parameters); see 
\cite[Chap.~15]{Lusztig03}. In Lemma~\ref{typedn}, we shall show that
this also holds if $r=0$.

Recall that we have a natural parametrisation of 
$\Irr(W)$ by the set $\cP_m$ of all pairs of partitions $(\alpha,\beta)$ 
such that $|\alpha|+|\beta|=m$. Let us write 
\[ \Irr(W)=\{E^{(\alpha,\beta)} \mid (\alpha,\beta)\in \cP_m\}.\]
By \cite[Prop.~22.14]{Lusztig03}, $f_E$ is a power of $2$ for all 
$E\in \Irr(W)$. So let us consider the subring $R=A_2 \subseteq {\C}(v)$ 
defined in Example~\ref{exp1}(c), where $p=2$. The ring $R$ is local, 
with maximal ideal $2R$, residue field $k=R/2R={\F}_2(v)$ and field of 
fractions $F={\Q}(v)$. As explained in (\ref{remcent}), the canonical
map $R \rightarrow k$ induces a decomposition map 
\[ d_{r,2} \colon R_0(H_F)\rightarrow R_0(H_k).\]
Note that $H_F$ and $H_k$ are split and $H_F$ is semisimple; see 
Dipper--James--Murphy \cite[\S 4, \S 6]{DJM}. Let $D_{r,2}$ be the 
corresponding decomposition matrix. Since $R$ is $\bc$-adapted for every  
two-sided cell $\bc$ of $W$, we have 
\begin{equation*}
D_{r,2}=\Big([E:[\Gamma_V]]\Big)_{E\in \Irr(W),V\in \Irr(H_k)};
\qquad \mbox{see (\ref{remcent})(3)}. 
\end{equation*}
Note that, by (\ref{remcent})(2), the columns of $D_{r,2}$ are linearly
independent. In fact, by Dipper--James--Murphy \cite[6.5]{DJM}, something
stronger holds: the rows and columns may be ordered so that $D_{r,2}$ 
has a lower triangular shape, with $1$ on the diagonal.

The discussion in (\ref{remcent}) also showed that, given 
any left cell $\Gamma$, there exists a unique $V \in \Irr(H_k)$ such 
that $[\Gamma] \cong [\Gamma_V]$. 

Now consider the constructible representations of $W$ with respect to
a weight function $L$ as above. The set $\Con(W)$ and the corresponding
``decomposition matrix'' matrix $\bD_r$ are determined by Lusztig 
\cite[22.24, 22.25 and 22.26]{Lusztig03} in a purely combinatorial way. 
We shall need the following result.
\end{abschnitt}

\begin{prop} \label{linind} Let $W$ be of type $B_m$ ($m \geq 2$) and
$L$ be a weight function as in (\ref{typeB}). Then there exist rational 
numbers $n_P$ ($P\in \Con(W)$) such that the following identity holds:
\[\dim E=\sum_{P\in \Con(W)} n_P\, [E:P] \qquad \mbox{for any 
$E\in \Irr(W)$}.\]
\end{prop}

The following proof is due to B.~Leclerc and reproduced here with his kind
permission. It requires some standard results about highest weight modules 
for the Lie algebra $\mathfrak{sl}_{n+1}(\C)$; see \cite{HyLe} and 
the references there. (This proof replaces my earlier proof which also 
used results from \cite{HyLe} but which was much less elementary.) 

\begin{proof} First note that the above identity simply means that the
regular representation of $W$ can be written as a rational linear 
combination of constructible representations. To be precise, we have to work
here in the appropriate Grothendieck group and identify a ${\C}[W]$-module
$\tilde{E}$ with 
\[ \sum_{E\in \Irr(W)} [E:\tilde{E}]\, E \in {\Q}[\Irr(W)].\]
To prove the above assertion, we place 
ourselves in the setting of Leclerc--Miyachi \cite{HyLe}, where 
the constructible representations of $W$ are interpreted in terms of 
canonical bases. Choose large positive integers $k,n$. (For example, 
any $k \geq m$ and $n\geq m-1+k+r$ will do; see \cite[4.1.1]{HyLe}.)
We consider the simple Lie algebra $\fg=\mathfrak{sl}_{n+1}(\C)$,
with Chevalley generators $e_j,h_j,f_j$ ($1\leq j \leq n$). Let 
$V(\Lambda)$  be the highest weight module of $\fg$ with highest weight 
$\Lambda=\Lambda_k+\Lambda_{k+r}$, where $\Lambda_1,\ldots,\Lambda_n$
are the fundamental weights\footnote{Note that, in \cite{HyLe}, the
quantized versions of these modules are considered. For our purposes,
it is enough to work with the specialisations at $1$.}. There is a canonical 
embedding of $\fg$-modules $\Phi \colon V(\Lambda) \rightarrow F(\Lambda)$ 
where $F(\Lambda)=V(\Lambda_{k+r}) \otimes V(\Lambda_k)$. The module 
$F(\Lambda)$ has a standard basis 
\[ \cS(\Lambda)=\{u_S\mid S\in \mbox{Sy}(n,k,r)\}\]
where $\mbox{Sy}(n,k,r)$ is the set of all symbols of the form 
\begin{equation*}
S=\left(\begin{array}{c} 1 \leq \beta_1 <\cdots < \beta_{k+r} \leq n+1\\
1\leq \gamma_1<\cdots <\gamma_k\leq n+1 \end{array}\right);\tag{1}
\end{equation*}
see \cite[2.2]{HyLe}. The action of the Chevalley generator $f_j$ on 
a basis element $u_{S}$ is explicitly given by the formulas in 
\cite[2.1]{HyLe}. One notices that these formulas constitute a refinement 
of the ``branching rule'' (see \cite[6.1.9]{ourbuch}), that is, we have  
\begin{equation*}
\mbox{ind}(u_S)=f_1(u_S)+f_2(u_S)+\cdots +f_n(u_S),\tag{2}
\end{equation*}
where $\mbox{ind}(u_S)$ denotes the sum of all $u_{S'}$ where $S'$ is a
symbol as in (1) which can be obtained from $S$ by increasing exactly one 
entry by $1$. (This means: we may replace some $\beta_i$ by $\beta_i+1$ 
if the new list of coefficients still satisfies the conditions in~(1); or 
we may replace some $\gamma_i$ by $\gamma_i+1$ if the new list of
coefficients still satisfies the conditions in~(1).) 

Let us now turn to the module $V(\Lambda)$. There is a well-defined symbol 
$S_0$ such that $u_{S_0}=\Phi(v_{\Lambda})$ where $v_\Lambda$ is a highest 
weight vector in $V(\Lambda)$; see \cite[2.2]{HyLe}. Thus, we have
\[ \Phi(V(\Lambda))=\langle f_{i_1}\circ \cdots \circ f_{i_s}(u_{S_0})\mid 
s\geq 0, 1 \leq i_1,\ldots,i_s \leq n\rangle_{\C}\subseteq F(\Lambda).\]
A symbol $S$ as in (1) is called {\em standard} if $\beta_i\leq \gamma_i$ 
for $1\leq i\leq k$. Let $\mbox{SSy}(n,k,r)$ be the set of all standard 
symbols in $\mbox{Sy}(n,k,r)$. By \cite[2.3]{HyLe}, this set will label
the {\em canonical basis} (or {\em lower global basis}) $\cB(\Lambda)$ of 
$V(\Lambda)$; see, for example, Jantzen \cite{Ja} for the general theory 
of canonical bases. In the present situation, we can just accept the 
conditions in \cite[2.4]{HyLe} as a definition, and then 
\cite[Theorem~3]{HyLe} provides an explicit construction of $\cB(\Lambda)=
\{b_S \mid S \in \mbox{SSy}(n,k,r)\}$. That construction proceeds by 
induction on the ``principal degree'' $d(S)$ of a symbol $S$, which is
defined as 
\begin{equation*}
d(S)=\sum_i \beta_i+\sum_j \gamma_j-\binom{k+1}{2}-\binom{k+r+1}{2}.\tag{3}
\end{equation*}
Given $d\geq 0$, let $\mbox{Sy}(n,k,r,d)$ be the set of symbols $S$ as
in (1) such that $d(S)=d$. We also write 
\[\mbox{SSy}(n,k,r,d):=\mbox{Sy}(n,k,r,d) \cap\mbox{SSy}(n,k,r).\]
The Leclerc--Miyachi construction shows that, for any $S\in 
\mbox{SSy}(n,k,r,d)$, there exists an explicitly known subset
$\cC(S)\subseteq \mbox{Sy}(n,k,r,d)$ such that
\begin{equation*}
\Phi(b_S)=\sum_{\Sigma\in \cC(S)} u_{\Sigma}.\tag{4}
\end{equation*}
Let us set $W_d=W(B_d)$ for $d=0,1,2,\ldots$. We bring now into 
play the representations of $W_d$. By \cite[4.1.1]{HyLe}, we may naturally 
label the irreducible representations of $W_d$ by the set of symbols in
$\mbox{Sy}(n,k,r,d)$. Thus, we write 
\begin{equation*}
\Irr(W_d)=\{E_S \mid S \in \mbox{Sy}(n,k,r,d)\}.
\end{equation*}
Identifying $u_S \leftrightarrow E_S$ for any $S$, we may just consider 
$\Irr(W_d)$ as a subset of $F(\Lambda)$ for any $d\geq 0$. With this 
convention, \cite[Theorem~10]{HyLe} states that 
\begin{equation*}
\Con(W_d)=\{\Phi(b_S) \mid S \in \mbox{SSy}(n,k,r,d)\} \qquad
\mbox{for $d=0,1,2, \ldots$}.\tag{5}
\end{equation*}
Now we have natural embeddings $\{1\}=W_0 \subset W_1\subset W_2 \subset 
\cdots$. Using our identification $u_S \leftrightarrow E_S$, we see that 
(2) does yield the ``branching rule'' for the induction of representations
from $W_d$ to $W_{d+1}$, that is, we have  
\[\Ind_{W_{d}}^{W_{d+1}}(E_S)=f_1(E_S)\oplus \cdots \oplus f_n(E_S)
\quad \mbox{for any $S\in \mbox{Sy}(n,k,r,d)$}.\]
Hence, starting with $\Irr(W_0)=\{{\bf 1}\}$, we obtain
\begin{equation*}
(f_1+\cdots +f_n)^d ({\bf 1})=\Ind_{W_0}^{W_d}({\bf 1})=\mbox{regular
representation of $W_d$}.\tag{6}
 \end{equation*}
Under the identification $u_S \leftrightarrow E_S$, the unit 
representation ${\bf 1}$ of $W_0$ corresponds to $u_{S_0}=\Phi(v_\Lambda)$. 
Now, the vector $(f_1+\cdots + f_n)^d(v_\Lambda)$ lies in the principal 
degree $d$ component of $V(\Lambda)$ and, hence, can be expressed as a 
rational linear combination of basis elements $b_S$ with $d(S)=d$. Applying 
$\Phi$ and using (5) and (6), we can now conclude that the regular 
representation of $W_d$ is a rational linear combination of $\Con(W_d)$, 
as desired. Taking $d=m$ (which is possible since $n,k$ were chosen 
sufficiently large) yields the desired assertion. 
\end{proof}

Note that, by Remark~\ref{con1}(b), the numbers $n_P$ are actually seen to 
be non-negative integers, once we have shown that Conjecture~\ref{conjC} 
holds. 

\begin{prop} \label{prop13} Let $W$ be of type $B_m$ ($m \geq 2$), with 
generators and relations given by the diagram in Example~\ref{expB1}. Let 
$L$ be a weight function such that $L(t)=b>0$ and $L(s_i)=a>0$ for $1\leq i 
\leq m-1$. Assume that $b=ar$ for some $r\geq 1$ and that (P1)--(P15) 
hold for $W,L$. Then $[\Gamma]$ is constructible for every left cell 
$\Gamma$ and all construcible representations arise in this way. Thus,
Conjecture~\ref{conjC} holds for $W,L$.
\end{prop}

\begin{proof} It is easily shown (see, for example, 
\cite[Remark~2.15]{mykl}) that the left cell representations and the 
constructible representations only depend on $r$. Thus, we may assume 
without loss of generality that $a=1$ and $b=r\geq 1$. Now recall the 
description of the decomposition matrix $D_{2,r}$ in (\ref{typeB}). The 
representations carried by the left cells of $W$ are completely 
determined by $D_{r,2}$.

Now consider the constructible representations of $W$. Let $P\in \Con(W)$.
By \cite[Lemma~22.2]{Lusztig03}, there exists some left cell $[\Gamma]$
such that $P\cong [\Gamma]$ as ${\C}[W]$-modules. So we have $P\cong 
[\Gamma]\cong [\Gamma_V]$ for a unique $V \in \Irr(H_k)$. Thus, there 
exists a subset $\cC'\subseteq \Irr(H_k)$ such that 
\begin{equation*}
\Con(W)=\{[\Gamma_V] \mid V \in \cC'\}.\tag{1}
\end{equation*}
Now, by (\ref{remcent}), the modules $[\Gamma_V]_R$ ($V \in \Irr(H_k)$)
are precisely the PIM's of $H_R$ (up to isomorphism). So we have an
isomorphism of $H_R$-modules
\[ H_R \cong \bigoplus_{V \in \Irr(H_k)} (\dim V)\,[\Gamma_V]_R\]
where $H_R$ is regarded as an $H_R$-module via left multiplication. 
Extending scalars from $R$ to $K$ and considering the multiplicities
of simple modules, we find
\begin{align*}
\dim E &=[E:{\C}[W]]=[E_v:H_K]=\sum_{V \in \Irr(H_k)} (\dim V)\,
[E_v:[\Gamma_V]_K]\\&=\sum_{V \in \Irr(H_k)} (\dim V)\, [E:[\Gamma]]
\quad \mbox{for any $E \in \Irr(W)$}.  
\end{align*}
Using now (1) and Proposition~\ref{linind}, we obtain the identity
\begin{equation*}
\sum_{V \in \Irr(H_k)} (\dim V)\,[E:[\Gamma_V]]=\sum_{V\in \cC'} 
n_V\, [E:[\Gamma_V]] \quad\mbox{for any $E\in \Irr(W)$},\tag{2}
\end{equation*}
where $n_V$ are certain rational numbers. By (\ref{typeB}), the coefficients 
$[E:[\Gamma_V]]$ occurring in the above identity are the entries of 
$D_{r,2}$. Now, by (\ref{remcent})(2), the columns of $D_{r,2}$ are 
linearly independent over $\Q$. Thus, we can compare coefficients in (2). 
In particular, every summand on the left hand side must also occur on the
right hand side and so $\cC'= \Irr(H_k)$. This shows that all
representations carried by the left cells are constructible, as desired.
\end{proof}

\begin{rem} \label{redu}  Let $W$ be of type $B_m$ and $L$ be a weight 
function as in Proposition~\ref{prop13}, that is, we have $b=ra$ for 
some $r\geq 1$. 

(a) The above result shows that the decomposition matrix $D_{r,2}$ in
(\ref{typeB}) coincides (up to a permutation of the columns) with the 
``decomposition matrix'' $\bD_r$ giving the expansion of the 
constructible representations in terms of the irreducible ones. Thus,
we obtain a new proof of \cite[Theorem~16]{HyLe}, assuming that 
(P1)--(P15) hold. The proof in [{\em loc. cit.}] uses the deep results
of Ariki \cite{Ar}.

(b) Table~II in \cite[p.~35]{Lusztig77} shows that $L$ arises (in the 
sense explained by Lusztig \cite[Chap.~0]{Lusztig03}) in the representation 
theory of finite classical groups. Hence, there is some hope that the 
geometric realization of $H$ in \cite[Chap.~27]{Lusztig03} will lead to a 
proof of (P1)--(P15) in this case.
\end{rem}

\begin{abschnitt} {\bf Type $D_m$.} \label{typedn} Let $W_1$ be a Coxeter
group of type $D_m$ ($m\geq 3$) with generators and relations given by the 
following diagram:
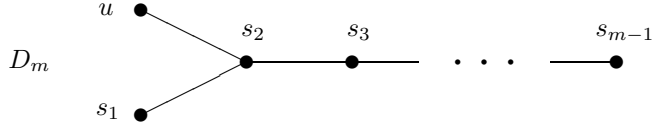
\begin{center}
\makeatletter
\begin{picture}(300,50)
\put( 0, 18){$D_m$}
\put( 50 , 0){\circle*{5}}
\put( 50 ,40){\circle*{5}}
\put( 50, 40){\line(2,-1){40}}
\put( 50,  0){\line(2,1){40}}
\put( 33,  0){$s_1$}
\put( 34, 38){$u$}
\put( 90 ,20){\circle*{5}}
\put(130 ,20){\circle*{5}}
\put(230 ,20){\circle*{5}}
\put( 90, 20){\line(1,0){40}}
\put(130, 20){\line(1,0){25}}
\put(170, 20){\circle*{2}}
\put(180, 20){\circle*{2}}
\put(190, 20){\circle*{2}}
\put(205, 20){\line(1,0){25}}
\put( 88 ,30){$s_2$}
\put(128, 30){$s_3$}
\put(222, 30){$s_{m{-}1}$}
\end{picture}
\makeatother
\end{center}
Let $\omega\colon W_1 \rightarrow W_1$ be the automorphism such that 
$\omega(u)=s_1$, $\omega(s_1)=u$ and $\omega(s_i)=s_i$ for $i>1$. Then the 
semidirect product $W=W_1 \rtimes \langle \omega\rangle$ can be naturally 
identified with a Coxeter group of type $B_m$, with generators $\{\omega,
s_1,s_2, \ldots,s_{m-1}\}$. Regarding $\omega$ as an element of $W$, we 
have the relation $\omega w= \omega(w)\omega$ for all $w\in W_1$. 
In particular, we have $u=\omega s_1 \omega$. 

Given any ${\C}[W_1]$-module $E$, we can define a new ${\C}[W_1]$-module 
structure on $E$ by composing the original action with the automorphism 
$\omega$. We denote that new ${\C}[W_1]$-module by $^\omega E$. 

Recall that $\Irr(W)=\{E^{(\alpha,\beta)} \mid (\alpha,\beta)\in \cP_m\}$
see (\ref{typeB}). For $(\alpha,\beta)\in \cP_m$, we denote by 
$E^{[\alpha,\beta]}$ the restriction of $E^{(\alpha,\beta)}$ to $W_1$. 
Then we have 
\begin{alignat*}{2}
E^{[\alpha,\beta]}&\cong E^{[\beta,\alpha]} \in \Irr(W_1) && \qquad \mbox{if
$\alpha\neq \beta$},\\
E^{[\alpha,\alpha]}&\cong E^{[\alpha,+]}\oplus E^{[\alpha,-]} && \qquad 
\mbox{if $\alpha=\beta$},
\end{alignat*}
where $E^{[\alpha,\pm ]}\in \Irr(W_1)$ and $E^{[\alpha,+]} \not\cong 
E^{[\alpha,-]}\cong {^\omega} E^{[\alpha,+]}$. This yields (see 
\cite[Chap.~5]{ourbuch} for more details): 
\[ \Irr(W_1)=\{E^{[\alpha,\beta]} \mid \alpha \neq \beta\}\cup
\{E^{[\alpha,\pm]} \mid 2|\alpha|=m\}.\]
Let $L\colon W \rightarrow \Z$ be the weight function such that 
\[ L(\omega)=0 \qquad \mbox{and}\qquad L(s_1)=L(s_2)=\cdots =L(s_{m-1})=1.\]
Let $L_1$ be the restriction of $L$ to $W_1$. Then $L_1$ is just the 
usual length function on $W_1$; see, for example, 
\cite[Lemma~1.4.12]{ourbuch}.

Let $H$ be the Iwahori--Hecke algebra associated with $W,L$. Let $H_1$ 
be the $A$-subspace of $H$ spanned by all $T_{w_1}$ with $w_1 \in W_1$. 
Then $H_1$ is nothing but the Iwahori--Hecke algebra associated with 
$W_1,L_1$. Note $T_\omega^2=T_1$ since $L(\omega)=0$. We have 
\[ T_{w_1}T_\omega=T_{w_1\omega} \quad \mbox{and}\quad T_{\omega}T_{w_1}=
T_{\omega(w_1)}T_\omega \quad \mbox{for any $w_1\in W_1$}.\]
In the following discussion, it will be understood that a left cell of 
$W_1$ is defined with respect to $L_1$ and a left cell of $W$ is defined 
with respect to $L$.
\end{abschnitt}

\begin{lem} \label{lemdn} In the setting of (\ref{typedn}), (P1)--(P15) 
hold for $W,L$. Let $\Gamma_1$ be a left cell of $W_1$. Then 
$\omega(\Gamma_1)$ also a left cell of $W_1$ and 
\[ \Gamma:=\Gamma_1 \cup \omega(\Gamma_1)\omega\]
is a left cell of $W$. Furthermore, we have isomorphisms of ${\C}[W]$-modules 
\[ \Res_{W_1}^W([\Gamma])\cong [\Gamma_1]\oplus [\omega(\Gamma_1)] \quad 
\mbox{and} \quad [\omega(\Gamma_1)]\cong {^\omega [\Gamma_1]}.\]
\end{lem}

\begin{proof} We consider the Kazhdan--Lusztig basis $\{c_w \mid w\in W\}$ 
of $H$ (defined with respect to $W,L$). Note that the construction of that
basis in \cite[Chap.~5]{Lusztig03} works for weight functions with
arbitrary integer values (even negative ones). We have $c_w=T_w+\sum_{y<w} 
p_{y,w} T_y$, where $p_{y,w} \in v^{-1}{\Z}[v^{-1}]$. Let $y,w\in W$ and 
write $y=y_1\omega^i$, $w=w_1 \omega^j$ where $y_1, w_1 \in W_1$ and $i,j\in 
\{0,1\}$. Then we have $p_{y,w}=0$ if $i\neq j$ and $p_{y,w}=p_{y_1,w_1}$ 
if $i=j$, where $p_{y_1,w_1}$ is the polynomial defined with respect to
$W_1,L_1$; see the remarks in \cite[\S 3]{Lu4} and \cite[\S 2]{my00}. Thus, 
$\{c_{w_1}\mid w_1\in W_1\}$ is the Kazhdan--Lusztig basis of $H_1$ (with 
respect to $L_1$). We have $c_\omega=T_\omega$ and 
\[ c_{w_1}c_\omega=c_{w_1\omega} \quad \mbox{and}\quad c_{\omega}c_{w_1}=
c_{\omega(w_1)}c_\omega \quad \mbox{for any $w_1\in W_1$}.\]
Furthermore, let $\ba(z)$ ($z\in W$) be defined with respect to $L$, and let
$\ba_1(z_1)$ ($z_1\in W_1$) be defined with respect to $L_1$. Then we have 
\[ \ba(w_1\omega)=\ba(w_1)=\ba_1(w) \qquad \mbox{for any $w_1\in W$}.\]
Now (P1)--(P15) hold for $W_1,L_1$; see \cite[Chap.~15]{Lusztig03}. The
above relations imply that (P1)--(P15) also hold for $W,L$; see 
\cite[\S 3]{Lu4} where this is worked out explicitly. The above relations
also show that 
\begin{equation*}
\cD=\cD_1=\{z\in W\mid \ba(z)=\Delta(z)\}\subseteq W_1.\tag{1}
\end{equation*}
Let $J$ be the $\Z$-algebra with basis $\{t_w\mid w\in W\}$ defined with
respect to $W,L$. Then we have $t_\omega^2=t_1$ and 
\begin{equation*}
t_{w_1}t_\omega=t_{w_1\omega} \quad \mbox{and}\quad t_{\omega}t_{w_1}=
t_{\omega(w_1)}t_\omega \quad \mbox{for all $w_1\in W_1$}.\tag{2}
\end{equation*}
Let $J_1$ be the $\Z$-submodule of $J$ spanned by all elements $t_{w_1}$
for $w_1 \in W_1$. Then $J_1$ is nothing but the $\Z$-algebra defined
with respect to $W_1,L_1$. Furthermore, the $A$-algebra homomorphism 
$\phi \colon H \rightarrow J_A$ restricts to the homomorphism $\phi_1 
\colon H_1 \rightarrow (J_1)_A$ defined with respect to $W_1,L_1$. 
These remarks allow us to switch back and forth between left cells of
$W_1$ and left cells of $W$. For example, the characterisation of 
$\sim_{\cL}$ in Remark~\ref{rem0} shows that 
\begin{align*}
&\mbox{every left cell of $W_1$ is contained in a left cell of $W$};\tag{3a}\\
&\mbox{if $\Gamma_1$ is a left cell of $W_1$, then so is $\omega(\Gamma_1)$}.
\tag{3b}
\end{align*}
Now let $\Gamma_1$ be a left cell of $W_1$. By (3b), we know that 
$\omega(\Gamma_1)$ is a left cell of $W_1$. By (3a), we have $\Gamma_1
\subseteq \Gamma$ where $\Gamma$ is a left cell of $W$. First we claim that 
$\Gamma \cap W_1=\Gamma_1$. To see this, let us fix an element $z\in 
\Gamma_1$. Now let $y\in \Gamma \cap W_1$. Since $z,y \in\Gamma$, we have 
$t_{y}t_{z^{-1}}\neq 0$ inside $J$, by Remark~\ref{rem0}. Since $z,y\in W_1$,
we also have $t_{y}t_{z^{-1}} \neq 0$ inside $J_1$ and so $y\in \Gamma_1$ 
(again using Remark~\ref{rem0}). Thus, the above claim is proved. Now let 
$y\in \Gamma$ and assume that $y \not\in W_1$. Let us write $y= y_1\omega$ 
where $y_1 \in W_1$. Since $y\in \Gamma$, we have $t_z t_{y^{-1}} \neq 0$ 
inside $J$.  Using (2), we also have $t_z t_{\omega(y_1)^{-1}} \neq 0$ 
inside $J_1$ and so $\omega(y_1) \in \Gamma_1$. Thus, we have shown that 
$\Gamma=\Gamma_1\cup \omega(\Gamma_1)\omega$, as desired. 

Finally, consider the statement concerning the left cell representations.
Let $J_{\C}^\Gamma:=\langle t_y\mid y \in\Gamma\rangle_{\C}\subseteq 
J_{\C}$.  By Lemma~\ref{lem30}, we have $J_{\C}^\Gamma=J_{\C}t_d$ where
$\cD\cap \Gamma=\{d\}$. By \cite[Lemma~21.2]{Lusztig03}, we have 
\[ J_{\C}^\Gamma \cong [\Gamma]_{\spadesuit} \qquad \mbox{(as 
${\C}[W]$-modules)}.\]
A similar statement holds for $\Gamma_1\subseteq W_1$. Hence, we can 
translate the desired assertion about the restriction of $[\Gamma]$ from
$W$ to $W_1$ to a statement about the restriction of $J_{\C}^\Gamma$ from
$J_{\C}$ to $(J_1)_{\C}$. Now, the decomposition $\Gamma=\Gamma_1\cup 
\omega(\Gamma_1)\omega$ yields that 
\[ (J_{\C}^\Gamma)\cong (J_1)_{\C}^{\Gamma_1}\oplus 
(J_1)_{\C}^{\omega(\Gamma_1)} \qquad \mbox{(as $\C$-vectorspaces)}.\]
But then (2) shows that this is an isomorphism of $(J_1)_{\C}$-modules.
Furthermore, for any $w\in W_1$, the action of any $t_w$ on 
$(J_1)_{\C}^{\omega(\Gamma_1)}$ will be the same as the action of
$t_{\omega(w)}=t_\omega t_wt_\omega$ on $(J_1)_{\C}^{\Gamma_1}$, as required.
\end{proof}

\begin{prop} \label{conD} Let $W_1$ be of type $D_m$, as in (\ref{typedn}), 
and let $L_1$ be the weight function on $W_1$ given by the length. By 
\cite[Chap.~15]{Lusztig03}, (P1)--(P15) hold for $W,L_1$. Then $[\Gamma_1]$ 
is constructible for every left cell $[\Gamma_1]$, and all constructible
representations arise in this way. Thus, Conjecture~\ref{conjC} holds for
$W_1,L_1$.
\end{prop}

\begin{proof} By \cite[Lemma~22.2]{Lusztig03}, we already know that
every constructible representation is of the form $[\Gamma_1]$ for some
left cell $\Gamma_1$. Conversely, let $\Gamma_1$ be a left cell of
$W_1$. We must show that $[\Gamma_1]$ is constructible. For this purpose, 
we place ourselves in the setting of Lemma~\ref{lemdn}.

Assume first that $[\Gamma_1]$ has a simple component $E_1$ with $f_{E_1}=1$.
Then Corollary~\ref{cor12a} shows that $[\Gamma_1]$ is irreducible and 
constructible. Now assume that $[\Gamma_1]$ has no simple component $E_1$ 
with $f_{E_1}=1$. Using the formulas for $f_E$ in \cite[22.14]{Lusztig03}, 
we see that all simple components of $[\Gamma_1]$ must be of the form 
$E^{[\alpha,\beta]}$ with $\alpha\neq \beta$. Thus, by Lemma~\ref{lemdn},
we have $[\omega(\Gamma_1)] \cong {^\omega [\Gamma_1]} \cong [\Gamma_1]$ and 
so 
\[ \Res_{W_1}^W([\Gamma])\cong 2[\Gamma_1].\]
Now, the argument in the proof of Proposition~\ref{prop13} also works
in the present situation, where $L(\omega)=0$ and $L(s_i)=1$ for all $i$.
Hence, $[\Gamma]$ is a constructible representation of $W$. Now
\cite[22.26]{Lusztig03} (see also Leclerc--Miyachi \cite[Theorem~13]{HyLe}) 
shows that $\Res_{W_1}^W([\Gamma])$ is twice a constructible representation 
of $W_1$. Hence $[\Gamma_1]$ is constructible as desired.
\end{proof}

\section{Exceptional types} \label{exccells}
We are now going to indicate proofs of Conjecture~\ref{conjC} for groups
of exceptional type. If $W$ is of type $E_6$, $E_7$, $E_8$ and $F_4$
with equal parameters, Lusztig's proof \cite{Lusztig86} requires some
sophisticated results from the representation theory of reductive groups 
over finite field, see \cite[Chap.~12]{LuBook}. We will show here that 
these arguments can be replaced by more elementary ones, involving some 
explicit computations with the character tables in \cite{ourbuch}. 

\begin{abschnitt} {\bf Type $I_2(m)$.} \label{celli2m} Let $W$ be of 
type $I_2(m)$ ($m\geq 3$), that is, we have $W=\langle s_1,s_2\rangle$ where 
$s_1^2=s_2^2=(s_1s_2)^m=1$. A weight function $L$ is specified by $L(s_1)=
a>0$ and $L(s_2)=b>0$ where $a=b$ if $m$ is odd. 

The irreducible representations of $W$ are given as follows. We have
the unit representation $1_W$ and $\mbox{sgn}$. If $m$ is even, there
are two further $1$-dimensional representations, which we denote by 
$\mbox{sgn}_1$ and $\mbox{sgn}_2$. They are characterised by the condition
that $s_1$ acts as $-1$ in $\mbox{sgn}_2$ and $s_2$ acts as $-1$ in 
$\mbox{sgn}$. All other irreducible representations have dimension~$2$; see 
\cite[5.3.4]{ourbuch} for an explicit description of these representations.
We denote by $\tau$ the direct sum of all the $2$-dimensional representations.

For any $k \geq 0$, write $1_k=s_1s_2s_1\cdots $ ($k$ factors) and
$2_k=s_2s_1s_2 \cdots$ ($k$ factors). Then the following hold:
\begin{itemize}
\item[(i)] If $m$ is odd, then the left cells are
\[ \{1_0\},\quad \{2_m\}, \quad \{2_1,1_2,2_3,\ldots,1_{m-1}\}, 
\quad \{1_1,2_2,1_3, \ldots,2_{m-1}\}.\]
The left cell representations are $1_W$, $\mbox{sgn}$, $\tau$, $\tau$, 
respectively.
\item[(ii)] If $m$ is even and $a=b$,  then the left cells are
\[ \{1_0\},\quad \{2_m\},\quad \{2_1,1_2,2_3,\ldots,2_{m-1}\}, 
\quad \{1_1,2_2,1_3, \ldots,1_{m-1}\}.\]
The representations are $1_W$, $\mbox{sgn}$, $\mbox{sgn}_1\oplus \tau$, 
$\mbox{sgn}_2\oplus \tau$, respectively.
\item[(iii)] If $m$ is even and $b>a$,  then the left cells are
\begin{gather*}
\{1_0\},\qquad \{1_1\},\qquad \{2_{m-1}\},\qquad \{2_m\},\\
\{2_1,1_2,2_3,\ldots,2_{m-2}\}, \quad \{2_2,1_3,2_4,\ldots,1_{m-1}\}.
\end{gather*}
The representations are $1_W$, $\mbox{sgn}_2$, $\mbox{sgn}_1$, $\mbox{sgn}$, 
$\tau$, $\tau$, respectively. 
\end{itemize}
The left cells in all of the above cases are determined in 
\cite[Chap.~8]{Lusztig03} (see also \cite[Exc.~11.4]{ourbuch} for
the case $a\neq b$). These computations do not require any of the
conditions (P1)--(P15). The representations carried by the left 
cells are easily determined by an explicit computation; see \cite[\S 6]{my02}
for case~(iii). The constructible representations are listed in
\cite[\S 12]{Lusztig82b} ($a=b$) and \cite[\S 6]{my02} ($b>a$). In each 
case, the representations carried by the left cells are precisely the 
constructible ones. Thus, Conjecture~\ref{conjC} holds for $W$ and 
any weight function $L$.
\end{abschnitt}

\begin{exmp} \label{cellh3} Let $W$ be of type $H_3$, with generators 
and relations given by the following diagram:
\begin{center}
\begin{picture}(220,20)
\put( 10, 5){$H_3$}
\put( 61,13){$s_1$}
\put( 91,13){$s_2$}
\put(121,13){$s_3$}
\put(78,9){$\scriptstyle{5}$}
\put( 65, 5){\circle*{5}}
\put( 95, 5){\circle*{5}}
\put(125, 5){\circle*{5}}
\put( 65, 5){\line(1,0){30}}
\put( 95, 5){\line(1,0){30}}
\end{picture}
\end{center}
All generators of $W$ are conjugate so every weight function on $W$ is of 
the form $L=al$ for some $a>0$. We have $|W|=120$. Using {\sf CHEVIE}
\cite{chevie}, one can explicitly compute the basis $\{c_w\}$ and all
polynomials $h_{x,y,z}$. By inspection, one sees that 
\[p_{yw}\in {\N}[v^{-1}]\qquad\mbox{and}\qquad h_{x,y,z}\in {\N}[v,v^{-1}]\]
for all $x,y,z,w\in W$. Thus, the arguments in \cite[Chap.~15]{Lusztig03}
show that (P1)--(P15) hold for $W,L$.

By \cite[Lemma~22.2]{Lusztig03}, we already know that each constructible
representation is of the form $[\Gamma]$ where $\Gamma$ is a left cell 
of $W$. The constructible representations are listed in 
\cite[\S 12]{Lusztig82b}. This yields the partition of $\Irr(W)$ into 
families. Using the tables for the $\Jind$-induction in 
\cite[Table~D.1]{ourbuch}, it is readily checked that there are no cuspidal 
families. We can now use the argument in (\ref{cusp}) to conclude that 
all left cell representations of $W$ are constructible. (Of course, having
computed all $p_{y,w}$ and all $h_{x,y,z}$, one could also directly compute
the left cells and the corresponding representations.) Thus, 
Conjecture~\ref{conjC} holds for $W$.
\end{exmp}

\begin{abschnitt} {\bf Type $H_4$.} \label{cellh4} Let $W$ be of type $H_4$,
with generators and relations given by the following diagram:
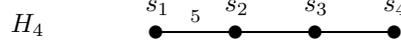
\begin{center}
\begin{picture}(220,20)
\put( 10, 5){$H_4$}
\put( 61,13){$s_1$}
\put( 91,13){$s_2$}
\put(121,13){$s_3$}
\put(151,13){$s_4$}
\put(78,9){$\scriptstyle{5}$}
\put( 65, 5){\circle*{5}}
\put( 95, 5){\circle*{5}}
\put(125, 5){\circle*{5}}
\put(155, 5){\circle*{5}}
\put( 65, 5){\line(1,0){30}}
\put( 95, 5){\line(1,0){30}}
\put(125, 5){\line(1,0){30}}
\end{picture}
\end{center}
All generators of $W$ are conjugate so every weight function on $W$ is of 
the form $L=al$ for some $a>0$. The constructible representations are 
determined by Alvis--Lusztig \cite{AlLu}. On the other hand, the left 
cells and the corresponding representations have been explicitly computed 
by Alvis \cite{Al}. By inspection, one sees that the representations carried 
by the left cells are precisely the constructible ones; see 
\cite[Prop.~3.5]{Al}. Thus, Conjecture~\ref{conjC} holds for $W$.
\end{abschnitt}

\begin{abschnitt} {\bf Type $F_4$.} \label{cellf4} 
Let $W$ be of type $F_4$, with generators and relations given by the 
following diagram:
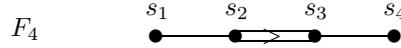
\begin{center}
\begin{picture}(220,20)
\put( 10, 5){$F_4$}
\put( 61,13){$s_1$}
\put( 91,13){$s_2$}
\put(121,13){$s_3$}
\put(151,13){$s_4$}
\put( 65, 5){\circle*{5}}
\put( 95, 5){\circle*{5}}
\put(125, 5){\circle*{5}}
\put(155, 5){\circle*{5}}
\put(105,2.5){$>$}
\put( 65, 5){\line(1,0){30}}
\put( 95, 7){\line(1,0){30}}
\put( 95, 3){\line(1,0){30}}
\put(125, 5){\line(1,0){30}}
\end{picture}
\end{center}
A weight function $L$ on $W$ is specified by two positive integers 
$a:=L(s_1)=L(s_2)>0$ and $b:=L(s_3)=L(s_4)>0$; we shall write $L=L_{a,b}$. 
In \cite{mykl}, the partition of $W$ into left cells has been determined 
for all values of $a,b$. Let $L=L_{a,b}$ and $L'=L_{a',b'}$ be two weight 
functions on $W$ such that $b\geq a>0$ and $b'\geq a'>0$. Then $L,L'$ define 
the same partition of $W$ into left cells if and only if $L,L'\in \cL_i$ for 
$i\in \{0,1,2,3\}$, where $\cL_i$ are defined as follows:
\begin{alignat*}{3}
\cL_0 &= \{ (c,c,c,c) && \mid c>0\}, \\
\cL_1 &= \{ (c,c,2c,2c) &&\mid c>0\},\\
\cL_2 &= \{ (c,c,d,d) &&\mid 2c>d>c>0\},\\
\cL_3 &= \{ (c,c,d,d) &&\mid d>2c>0\}.
\end{alignat*}
Furthermore, if $L,L'\in \cL_i$, then the left cells give rise to the
same representations of $W$. On the other hand, the constructible 
representations are determined in \cite[22.27]{Lusztig03}\footnote{As
pointed out in \cite[Remark~4.10]{mykl}, there is an error in the list
of constructible representations in \cite[22.27, Case~1]{Lusztig03},
where $b=2a>0$: the representations denoted $1_3\oplus 2_1$ and $1_2
\oplus 2_2$ have to be removed from that list.}. By inspection, one sees 
that the left cell representations are precisely the constructible ones. 
Thus, Conjecture~\ref{conjC} holds for $W$ and all weight functions $L$.

Note that, by Table~II in \cite[p.~35]{Lusztig77}, only the weight functions
$L$ such that $b=a$, $b=2a$ or $b=4a$ arise in the representation theory of
finite reductive groups.

For the case of equal parameters, let us indicate an argument which does
not require the explicit computation of all left cells. This will also 
be a model for the discussion of groups of type $E_6$, $E_7$ and $E_8$. So
let us assume that $L=al$ for some $a>0$. By \cite[Chap.~15]{Lusztig03}, the 
properties (P1)--(P15) hold for $W,L$. Hence, by \cite[Lemma~22.2]{Lusztig03},
we already know that every constructible representation is carried by a 
left cell of $W$. To prove the converse, it will now be enough to consider 
only those left cells which belong to a cuspidal family of $\Irr(W)$; see 
(\ref{cusp}).

By \cite[8.1]{LuBook}, there is a unique  cuspidal family $\cF_0$ of 
$\Irr(W)$, the one containing the representation $12_1$, where we use
the notation of \cite[4.10]{LuBook} or \cite[Table~C.3]{ourbuch}. To deal 
with this family, we consider the parabolic subgroup $W_I$ of type $B_3$, 
where $I=S\setminus \{s_4\}$. By Proposition~\ref{prop13}, the left cell 
representations of $W_I$ are precisely the constructible ones, and these 
are explicitly given in \cite[22.24]{Lusztig03}.

\begin{table}[htbp] 
\caption{The cuspidal family in type $F_4$} \label{cuspf4}
\begin{center}
$\renewcommand{\arraystretch}{1.1}
\begin{array}{cc|ccccc} \hline E & f_E &P_1 &P_2&P_3&P_4
&P_5\\ \hline
1_2 & 8 & 1 & . & . & . & . \\ 
1_3 & 8 & . & 1 & . & . & . \\
4_1 & 8 & . & . & 1 & . & . \\
4_3 & 4 & 1 & . & . & 1 & . \\
4_4 & 4 & . & 1 & . & . & 1 \\
6_1 & 3 & . & . & . & 1 & 1 \\
6_2 &12 & 1 & 1 & 1 & . & . \\
9_2 & 8 & 2 & . & 1 & 1 & . \\
9_3 & 8 & . & 2 & 1 & . & 1 \\
12_1&24 & 1 & 1 & 1 & 1 & 1 \\
16_1& 4 & 1 & 1 & 2 & 1 & 1 \\
\hline\end{array}$
\end{center}
\end{table}

Now let $\Gamma$ be a left cell of $W$ belonging to $\cF_0$.
Then, as in the proof of Lemma~\ref{lem42}, there exists a
unique left cell $\Gamma'$ of $W_I$ such that $\Gamma \subseteq
X_I\Gamma'$. Then $[\Gamma]$ is a direct summand of 
$\Ind_I^S([\Gamma'])$. Let us write  
\[ \Ind_I^S([\Gamma'])\cong \bigoplus_{E \in \cF_0} m_E E \;\oplus 
\mbox{ combination of $E'\in\Irr(W)\setminus \{\cF_0\}$},\]
where $m_E$ are non-negative integers. Thus, we can conclude that 
\[ [\Gamma]\cong \bigoplus_{E \in \cF_0} n_E E \qquad \mbox{where $0\leq 
n_E \leq m_E$ for all $E\in \cF_0$}.\]
Now the idea is to look for arithmetical conditions on the numbers $n_E$ 
so that the only remaining possibilities satisfying these conditions
correspond to the expansion of the constructible representations in $\cF_0$.
One such condition is given by Lemma~\ref{lem12}: the numbers $n_E$ must 
satisfy 
\begin{equation*}
\sum_{E\in \cF_0} \frac{1}{f_E}\, n_E =1.\tag{$\diamondsuit$}
\end{equation*}
Inducing all constructible representations from $W_I$ to $W$, we can
explicitly determine (using {\sf CHEVIE} \cite{chevie}) all possible
non-zero vectors $(m_E)_{E\in \cF_0}$ as above. They are given by the 
columns labelled by $P_1,P_2,P_3,P_4,P_5$ in Table~\ref{cuspf4}.
Now, if the vector $(m_E)_{E\in \cF_0}$ is given by one of the above
columns, we see that $\sum_E m_E/f_E=1$. Thus, ($\diamondsuit$) shows 
that we must have $n_E=m_E$ for all $E$. So $[\Gamma]$ is given by one 
of the above five columns. Comparison with the table of constructible 
respresentations in \cite[p.~223]{Lusztig82b} shows that $[\Gamma]$ is 
constructible. 
\end{abschnitt}

\begin{abschnitt} {\bf Type $E_6$.} \label{celle6} 
Let $W$ be of type $E_6$, with generators and relations given by the 
following diagram:
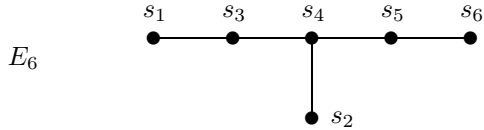
\begin{center}
\begin{picture}(250,50)
\put( 10,25){$E_6$}
\put(132, 3){$s_2$}
\put( 61,43){$s_1$}
\put( 91,43){$s_3$}
\put(121,43){$s_4$}
\put(151,43){$s_5$}
\put(181,43){$s_6$}
\put(125, 5){\circle*{5}}
\put( 65,35){\circle*{5}}
\put( 95,35){\circle*{5}}
\put(125,35){\circle*{5}}
\put(155,35){\circle*{5}}
\put(185,35){\circle*{5}}
\put(125,35){\line(0,-1){30}}
\put( 65,35){\line(1,0){30}}
\put( 95,35){\line(1,0){30}}
\put(125,35){\line(1,0){30}}
\put(155,35){\line(1,0){30}}
\end{picture}
\end{center}
All generators of $W$ are conjugate so every weight function on $W$ is of 
the form $L=al$ for some $a>0$. By \cite[Chap.~15]{Lusztig03}, the 
properties (P1)--(P15) hold for $W,L$. Hence, by \cite[Lemma~22.2]{Lusztig03},
we already know that every constructible representation is carried by a 
left cell of $W$. To prove the converse, it will now be enough to consider 
only those left cells which belong to a cuspidal family of $\Irr(W)$; see 
(\ref{cusp}). We shall argue as in (\ref{cellf4}).

By \cite[8.1]{LuBook}, there is a unique  cuspidal family $\cF_0$ of 
$\Irr(W)$, the one containing the representation $80_s$, where we use 
the notation of \cite[4.11]{LuBook} or \cite[Table~C.4]{ourbuch}. To deal 
with this family, we consider the parabolic subgroup $W_I$ of type $D_5$, 
where $I=S\setminus \{s_6\}$. By Proposition~\ref{conD}, the left cell 
representations of $W_I$ are precisely the constructible ones, and these 
are explicitly given in \cite[22.26]{Lusztig03}. The possibilities for 
all non-zero vectors $(m_E)_{E\in \cF_0}$ obtained (as above) by inducing 
all constructible representations of $W_I$ are given by the following table:
\[ \renewcommand{\arraystretch}{1.1}
\begin{array}{cc|ccc} \hline E & f_E &P_1 &P_2&P_3\\ \hline
10_s & 3 & 1 & . & . \\
20_s & 6 & . & 1 & . \\
60_s & 2 & 1 & . & 1 \\
80_s & 6 & 1 & 1 & 1 \\
90_s & 3 & . & 2 & 1\\ \hline\end{array}\]
Checking ($\diamondsuit$), we find that $(n_E)_{E\in \cF_0}$ must be equal 
to one of the vectors $(m_E)_{E\in \cF_0}$. So $[\Gamma]$ is given by one 
of the above three columns. Comparison with the table of constructible 
respresentations in \cite[p.~223]{Lusztig82b} shows that $[\Gamma]$ is
constructible. Thus, Conjecture~\ref{celle6} holds for $W$.
\end{abschnitt}

\begin{abschnitt} {\bf Type $E_7$.} \label{celle7} 
Let $W$ be of type $E_7$, with generators and relations given by the 
following diagram:
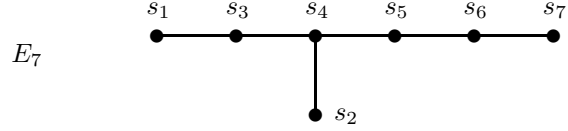
\begin{center}
\begin{picture}(270,50)
\put( 10,25){$E_7$}
\put(132, 3){$s_2$}
\put( 61,43){$s_1$}
\put( 91,43){$s_3$}
\put(121,43){$s_4$}
\put(151,43){$s_5$}
\put(181,43){$s_6$}
\put(211,43){$s_7$}
\put(125, 5){\circle*{5}}
\put( 65,35){\circle*{5}}
\put( 95,35){\circle*{5}}
\put(125,35){\circle*{5}}
\put(155,35){\circle*{5}}
\put(185,35){\circle*{5}}
\put(215,35){\circle*{5}}
\put(125,35){\line(0,-1){30}}
\put( 65,35){\line(1,0){30}}
\put( 95,35){\line(1,0){30}}
\put(125,35){\line(1,0){30}}
\put(155,35){\line(1,0){30}}
\put(185,35){\line(1,0){30}}
\end{picture}
\end{center}
In order to prove Conjecture~\ref{conjC}, we argue as in (\ref{celle6}).
It remains to consider the cuspidal families. By \cite[8.1]{LuBook}, 
there is a unique  cuspidal family $\cF_0$ of $\Irr(W)$, the one 
containing the representation $512_a'$, where we use the notation of 
\cite[4.12]{LuBook} or \cite[Table~C.5]{ourbuch}. To deal with this family, 
we consider the parabolic subgroup $W_I$ of type $E_6$.  The possibilities
for all non-zero vectors $(m_E)_{E\in \cF_0}$ obtained by inducing all 
constructible representations from $W_I$ are given by the following table:
\[ \renewcommand{\arraystretch}{1.1}
\begin{array}{cc|ccc} \hline E & f_E & P_1&P_2 &P_3 
\\ \hline 512_a & 2  & 1 & 2 & 3 \\
512_a' & 2 & 1 & 2 & 3 \\ \hline\end{array}\]
Now let $\Gamma$ be a left cell belonging to $\cF_0$. Checking 
($\diamondsuit$), we find the following possibilities for the expansion 
of $[\Gamma]$: $512_a \oplus 512_a'$, $2\cdot 512_a$ or $2 \cdot 512_a'$. The
first is constructible by the table in \cite[p.~223]{Lusztig82b}. We must 
show that the other two possibilities do not occur. This can be seen as 
follows. By Theorem~\ref{fam}, the family $\cF_0$ corresponds to a two-sided 
cell $\bc$ of $W$. Since $f_E=2$ for all $E\in \Irr(\cF_0)$, we consider the 
ring $A_2$ as in Example~\ref{exp1}. In fact, since we do not yet know if
$H$ is split over the residue field of $A_2$ (which is just ${\F}_2(v)$), 
we work with a ring $R\supseteq A_2$ as in (\ref{split}). 
It is readily checked that $R$ is $\bc$-adapted. So Theorem~\ref{prop32} 
shows that $e_{\bc}= \phi_R^{-1}(t_\bc)$ is a primitive idempotent in the 
center of $H_R$. Consider the decomposition matrix $D_{\bc}$ of the block 
$H_Re_{\bc}$. That matrix has only two rows, corresponding to $512_a$ and 
$512_a'$. The columns correspond to the projective indecomposable modules in 
that block. Now we already know that $512_a \oplus 512_a'$ is constructible 
and, hence, is carried by some left cell contained in $\bc$. Hence $512_a 
\oplus 512_a'$ gives one column of $D_{\bc}$; see Corollary~\ref{lem34}.
If $[\Gamma] \cong 2\cdot 512_a$ or $2\cdot 512_a'$, then this would give 
another column of $D_{\bc}$. Since $512_a'\cong 512_a \otimes \mbox{sgn}$, 
we would conclude that actually both $2\cdot 512_a$ and $2\cdot 512_a'$ 
correspond to columns of $D_{\bc}$. Thus, we would deduce that $D_{\bc}$ 
has at least three columns and, hence, does not have full rank. Consequently,
$D_R$ (which is a block diagonal matrix where one block is given by $D_{\bc}$)
would not have full rank either, contradicting (\ref{remcent})(2).
This completes the proof of Conjecture~\ref{conjC} for $W$.
\end{abschnitt}

\begin{table}[htbp] 
\caption{The cuspidal family in type $E_8$} \label{cusp8}
\begin{center}
$\renewcommand{\arraystretch}{1.1}
\begin{array}{ccc|cccccccc|ccccccc}
\hline E & (x,\sigma)& f_E & P_1&P_2&P_3&P_4&P_5&P_6&P_7&P_8&
\multicolumn{6}{c}{\mbox{splitting $P_8$}}\\\hline 
70_y& 1,\lambda^3&
30&1& .& .& .& .& .& .& .    & . &  . &  . &  . &  . &  . \\
168_y&g_2',\varepsilon' &
8 & .& 1& .& .& .& .& .& .   & . &  . &  . &  . &  . &  . \\
420_y&g_5,1&
5  & .& .& 1& .& .& .& .& .  & . &  . &  . &  . &  . &  . \\
448_w&g_2,\varepsilon&
12 &1& .& .& 1& .& .& .& .   & . &  . &  . &  . &  . &  . \\
1134_y&g_3,\varepsilon&
6 & .& .& .& 1& 1& .& .& .  & . &  . &  . &  . &  . &  . \\
1344_w&g_4,1&
4&.& 1& 1& .& 1& .& .& .    & . &  . &  . &  . &  . &  . \\
1400_y&1,\nu'&
24 & 2& 1& .& .& .& 1& .&  .& . &  . &  . &  . &  . &  . \\
1680_y&1,\lambda^2&
20 & 3& .& .& 1& .& 1& .&  .& . &  . &  . &  . &  . &  . \\
2688_y&g_2',\varepsilon''&
8& .& .& .& .& .& 1& 1& 2   & . &  2 &  . &  2 &  . &  2 \\
3150_y&g_3,1&
6 & .& .& 1& 1& 1& .& 1& 2  & 2 &  . &  2 &  . &  . &  2 \\
4200_y&g_2',1&
8 & .& 2& 1& .& 1& 1& 1& 2  & . &  2 &  . &  2 &  . &  2 \\
4480_y&1,1&
120 & 1& 1& 1& 1& 1& 1& 1& 2& 2 &  . &  2 &  . &  . &  2 \\
4536_y&1,\nu&
24 & 3& 1& .& 1& .& 2& 1& 2 & 2 &  . &  2 &  . &  . &  2 \\
5670_y&1,\lambda^1&
30&3& 1& .& 2& 1& 2& 1& 2   & 2 &  . &  2 &  . &  . &  2 \\
2016_w&g_6,1&
6&.& .& 1& .& .& .& 1& 2    & 2 &  . &  . &  2 &  2 &  . \\
5600_w&g_2,r&
6 &2& 1& .& 2& 1& 2& 1& 2   & . &  2 &  2 &  . &  2 &  . \\
7168_w&g_2,1&
12 & 1& 1& 1& 1& 1& 2& 2& 4 & 2 &  2 &  2 &  2 &  4 &  . \\ \hline 
\multicolumn{17}{l}{\mbox{The notation in the first two columns
is taken from Lusztig \cite[p.~369]{LuBook}.}}\\
\multicolumn{17}{l}{\mbox{The second column gives the 
identification with Gyoja's tables \cite[\S 2.3]{Gy}.}}
\end{array}$
\end{center}
\end{table}

\begin{abschnitt} {\bf Type $E_8$.} \label{celle8} 
Let $W$ be of type $E_8$, with generators and relations given by the 
following diagram:
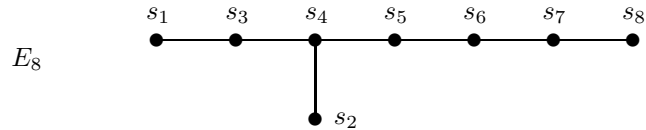
\begin{center}
\begin{picture}(280,50)
\put( 10,25){$E_8$}
\put(132, 3){$s_2$}
\put( 61,43){$s_1$}
\put( 91,43){$s_3$}
\put(121,43){$s_4$}
\put(151,43){$s_5$}
\put(181,43){$s_6$}
\put(211,43){$s_7$}
\put(241,43){$s_8$}
\put(125, 5){\circle*{5}}
\put( 65,35){\circle*{5}}
\put( 95,35){\circle*{5}}
\put(125,35){\circle*{5}}
\put(155,35){\circle*{5}}
\put(185,35){\circle*{5}}
\put(215,35){\circle*{5}}
\put(245,35){\circle*{5}}
\put(125,35){\line(0,-1){30}}
\put( 65,35){\line(1,0){30}}
\put( 95,35){\line(1,0){30}}
\put(125,35){\line(1,0){30}}
\put(155,35){\line(1,0){30}}
\put(185,35){\line(1,0){30}}
\put(215,35){\line(1,0){30}}
\end{picture}
\end{center}
In order to prove Conjecture~\ref{conjC}, we argue as in (\ref{celle6}).
It remains to consider the cuspidal families.
By \cite[8.1]{LuBook}, there is a unique  cuspidal family $\cF_0$ of 
$\Irr(W)$, the one containing the representation $4480_y$, where we use 
the notation of \cite[4.13]{LuBook} or \cite[Table~C.6]{ourbuch}. To deal 
with this family, we use the parabolic subgroup $W_I$ of type $E_7$. The 
possibilities for all non-zero vectors $(m_E)_{E\in \cF_0}$ cobtained by 
inducing all constructible representations from $W_I$ are given by the columns 
labelled $P_1,\ldots,P_8$ in Table~\ref{cusp8}. 

Checking ($\diamondsuit$) yields that $P_1,\ldots,P_7$ already give the 
expansion of a left cell representation in terms of the irreducible ones; 
furthermore, the table in \cite[p.~224]{Lusztig82b} shows that all these 
columns give constructible representations of $W$. 

Hence, it remains to consider the case where $\Gamma$ is a left cell 
such that $[\Gamma]$ is a direct summand of the column labelled $P_8$
in Table~\ref{cusp8}. We must show that the vector $(n_E)_{E\in \cF_0}$
giving the expansion of $[\Gamma]$ is obtained by dividing all coefficients 
in that column by~$2$. Now there are many more possibilities for 
$(n_E)_{E\in\cF_0}$ satisfying ($\diamondsuit$). So we have to look for
further arithmetical conditions on these numbers. Another such condition is
given by (\ref{center})(b), where we work over the ring $A_2$ as in 
Example~\ref{exp1}. The remaining possibilities for the vector 
$(n_E)_{E \in \cF_0}$ are listed in the $6$ rightmost columns of 
Table~\ref{cusp8}. (The central characters in (\ref{center}) can be 
computed using the programs contained in the file {\tt hecbloc.g} in the 
contributions directory on the {\sf CHEVIE} homepage \cite{chevie}.)

Now we can argue as follows.  Consider a ring $R\supseteq A_2$ as in 
(\ref{split}). Then, by (\ref{remcent})(2), the columns of the decomposition 
matrix of $H_R$ are linearly independent. Assume now, if possible, that 
$[\Gamma]$ is given by one of the $6$ rightmost columns in Table~\ref{cusp8}.
The $H_R$-module $[\Gamma]_R$ is projective (see Corollary~\ref{lem34}), hence 
it can be written as a direct sum of projective indecomposable $H_R$-modules.
(Note that $R$ is not $\bc$-adapted, where $\bc$ is the two-sided cell
containing $\Gamma$; so we don't know if $[\Gamma]_R$ is indecomposable.) 
Now we already know that there is a left cell $\Gamma_7$ whose
expansion in terms of irreducible representations is given by $P_7$, and 
$[\Gamma_7]_R$ is a projective $H_R$-module (Corollary~\ref{lem34}). Since 
$P_8=2P_7$, we conclude (using the linear independence of the columns of 
the decomposition matrix of $H_R$) that $[\Gamma]_R$ must be a sum of 
projective indecomposable $H_R$-modules which occur in the decomposition of 
$[\Gamma_7]_R$ as a sum of projective indecomposable $H_R$-modules. We 
don't know that decomposition of $[\Gamma_7]_R$ but we can just compute
all possible decompositions which satisfy the conditions in 
(\ref{center})(b). It turns out that the possibilities are precisely 
those given by the table with heading $(\Sym_5,\Sym_3 \times \Sym_2)$ 
($p=2$) of Gyoja \cite[p.~321]{Gy}. Thus, at least one of the $6$ rightmost
columns in Table~\ref{cusp8} should be expressible as a sum of the columns
in Gyoja's table. One easily checks that this impossible. This contradiction 
shows that Conjecture~\ref{conjC} holds for $W$.
\end{abschnitt}


\end{document}